\documentclass{article}[12pt]

\usepackage{amsmath}
\usepackage{amssymb}
\usepackage{amsfonts}
\usepackage{amsthm}
\usepackage[english]{babel}
\usepackage{graphicx}
\usepackage{color}
\usepackage{fancyhdr} 
\usepackage[left=1in,right=1in,top=1in,bottom=1in]{geometry}
\usepackage{amsmath,amscd}
\usepackage[all]{xy}
\usepackage{tikz}
\usetikzlibrary{shapes}
\usetikzlibrary{matrix}
\usepackage{mathtools}
\usepackage{caption}
\usepackage{hyperref}
\usepackage[numbers]{natbib}

\newcommand{\lrp}[1]{\left(#1\right)}
\newcommand{\lrb}[1]{\left[#1\right]}
\newcommand{\lrm}[1]{\left|#1\right|}
\newcommand{\lrc}[1]{\left\{#1\right\}}
\newcommand{\lra}[1]{\langle{#1}\rangle}

\newcommand{\eq}[2]{\begin{equation}\label{#2} \begin{split} #1  \end{split} \end{equation}}
\newcommand{\eqn}[1]{\begin{equation*} \begin{split} #1 \end{split} \end{equation*}}

\newcommand{\R}{\mathbb{R} }
\newcommand{\C}{\mathbb{C} }

\newcommand{\Z}{\mathbb{Z} }

\newcommand{\gv}{\mathbf{g} }
\newcommand{\vb}[1]{\mathbf{#1}}
\newcommand{\cA}{\mathcal{A} }
\newcommand{\cAp}{\mathcal{A}_{\mathrm{prin}} }
\newcommand{\cAps}[1]{\mathcal{A}_{\mathrm{prin},{#1}} }
\newcommand{\cX}{\mathcal{X} }
\newcommand{\ssO}{\mathcal{O} }

\newcommand{\mf}[1]{\mathfrak{#1} }

\newcommand{\matdd}[4]{\left(
        \begin{array}{c c}
            #1 & #2 \\ \noalign{\medskip}
            #3 & #4
        \end{array}
    \right)}

\DeclareMathOperator{\GL}{GL}
\DeclareMathOperator{\SL}{SL}
\DeclareMathOperator{\Mat}{Mat}

\DeclareMathOperator{\chara}{char}
\DeclareMathOperator{\gr}{gr}
\DeclareMathOperator{\diag}{diag}
 
\theoremstyle{plain}
\newtheorem{theorem}{Theorem}
\newtheorem{prop}[theorem]{Proposition}
\newtheorem{lemma}[theorem]{Lemma}

\newtheorem{cor}[theorem]{Corollary}

\theoremstyle{definition}

\newtheorem{example}[theorem]{Example}

\theoremstyle{remark}
\newtheorem{remark}[theorem]{Remark}

\newenvironment{problem}[1]{ \flushleft \textcolor{blue}{\normalsize {#1}}}

\newenvironment{subprob}[1]{ \flushleft \textcolor{blue}{\normalsize {#1}}}

\mathtoolsset{centercolon}
\begin{document}

\renewcommand{\thesubsection}{\arabic{section}.\alph{subsection}}
\renewcommand{\labelenumi}{(\arabic{enumi})}

\title{Fock-Goncharov conjecture and polyhedral cones for $U \subset \SL_n$ and base affine space $\SL_n/U$} 
\author{Timothy Magee\thanks{tmagee@math.utexas.edu, University of Texas at Austin}}
\date{}
\maketitle

{\abstract{I prove several conjectures of \cite{GHKK} on the cluster structure of $\SL_n$, 
which in particular imply the full Fock-Goncharov conjecture
for the open double Bruhat cell $\cA \subset \SL_n/U$, 
for $U \subset \SL_n$ a maximal unipotent subgroup. 
This endows the mirror cluster variety $\cX$ with a canonical
potential function $W$, 
and determines a canonical cone 
$W^T \geq 0  \subset \cX\lrp{\R^T}$ 
of the mirror tropical space, 
whose integer points parameterize a basis of
$H^0\lrp{\SL_n/U,\ssO_{\SL_n/U}}$, 
canonically determined by the open subset $\cA \subset \SL_n/U$.
Each choice of seed identifies $\cX\lrp{\R^T}$ with a real vector space, and
$W^T \geq 0$ with a system of linear equations with integer coefficients,
cutting out a polyhedral cone. 
We obtain in this way (generally) infinitely many parameterizations of the
canonical basis as integer points of a polyhedral cone. 
For the usual initial seed of the double Bruhat cell,
we recover the parameterizations of Berenstein-Kazhdan\cite{BKaz,BKaz2} and Berenstein-Zelevinsky\cite{BZ96} 
by integer points of the Gelfand-Tsetlin cone.
}}


\section{Introduction}
In this paper, I'll focus on two examples that are discussed in \cite{GHKK}.
First I'll look at a maximal unipotent subgroup $U$ of $G=\SL_n$, and then I'll consider the base affine space $G/U$.
I prove \cite{GHKK} Conjecture 11.11 and part 1 of Conjecture 0.19.
The notation and definitions used below come from \cite{GHKK}.

{\theorem{Every frozen variable of $G^{e,w_0}$ has an optimized seed.
Let $\mathring{U}$ be $G^{e,w_0}$ with the frozen variables $\Delta^{1,\dots,i}_{1,\dots,i}$ restricted to 1, 
and view $U$ as the partial compactification of $\mathring{U}$ afforded by allowing the remaining frozen variables to vanish.
Then the torus $T$ embedded in $U$ via the Lusztig map 
$\lrp{t_1,\cdots , t_{l\lrp{w_0}}} \mapsto E_{i_1}\lrp{t_1} \cdots E_{i_{l\lrp{w_0}}} \lrp{t_{l\lrp{w_0}}}$
for the lexographically minimal reduced word $w_0 = s_1 s_2 s_1 \cdots s_{n-1} s_{n-2}\cdots s_1$
is in $\mathring{U}$'s atlas of tori.
Each non-constant matrix entry is a cluster variable in some seed of $\mathring{U}$, 
and in the seed associated to $T$ the $\gv$-vectors for these non-constant matrix entries 
form a basis for the integer tropical points of the mirror
$\cX\lrp{\Z^T}$, which is identified with a lattice by the choice of seed.
The edges of the cone $W^T\geq 0 \subset \cX\lrp{\R^T}$ are the $\R_{\geq 0}$ span of these $\gv$-vectors,
so this cone is a full dimensional simplicial cone and the $\gv$-vectors are the non-identity generators of its monoid of integral points.
}}

{\theorem{\label{GTintro}
View $G/U$ as a partial compactification of $G^{e,w_0}$, 
where $G^{e,w_0}$ is embedded in $G/U$ via $g\mapsto g^T U$. 
In the seed determined by the lexicographically minimal reduced expression for $w_0$,\footnote{More precisely, the seed whose cluster variables are the chamber minors of the double pseudoline arrangement determined by $w_0 = s_1 s_2 s_1 s_3 s_2 s_1 \cdots$.}
$\Xi := \lrc{x \in \cX\lrp{\R^T} : W^T \geq 0}$ is the Gelfand-Tsetlin cone with the final coordinate restricted to 0.
The edges of $\Xi$ are the $\R_{\geq 0}$ span of the $\gv$-vectors of the minors $\Delta^{j_1, \dots, j_i}_{1,\dots,i}$.
These $\gv$-vectors are precisely the non-identity generators of the monoid $\Xi \cap \cX\lrp{\Z^T}$.}}

These theorems have some immediate consequences:

{\cor{The full Fock-Goncharov conjecture holds for $U \subset \SL_n$ and $\SL_n/U$.}}

{\cor{\label{dimVlambda}
Fix a weight $\lambda$, with $V_{\lambda}$ the associated irreducible representation of $\SL_n$.
$\lambda$ defines an ``affine subspace'' of $\cX\lrp{\R^T}$,
and its intersection with $\Xi$ is a bounded lattice ``polytope'' whose integer points parametrize a canonical basis for $V_{\lambda}$.
In particular, the number of integer points is $\dim V_\lambda$.}}

{\remark{The terms ``affine subspace'' and ``polytope'' are usually applied to linear spaces.
The real tropical space of an affine log Calabi-Yau is only piecewise linear, 
so it doesn't make sense to talk about a straight line in $\cX\lrp{\R^T}$.
Straight lines are replaced by ``broken lines''.
Each choice of seed identifies $\cX\lrp{\R^T}$ with a real vector space where the usual notions of ``affine subspace'' and ``polytope'' apply.
Under this identification, ``affine subspace'' and ``polytope'' take on their usual meaning. See \cite{GHKK} for details. 
}}

Corollary~\ref{dimVlambda} is of course already known for the Gelfand-Tsetlin cone.\cite{BZ96}
The proof of Corollary~\ref{dimVlambda} is, however, entirely different, at once conceptually
simpler, and much more broadly applicable. See the intro to \cite{GHKK} for
a discussion of the general picture.
 
Note that for both $U$ and $\SL_n/U$ the resulting canonical basis has a very natural set of functions.
For $U$ the canonical basis is generated by the most na\"ive thing possible-- the non-constant matrix entries.
For $\SL_{n}/U$, it is generated by the minors $\Delta^{j_1,\dots,j_i}_{1,\dots,i}$.  
%
%

{\remark{
In \cite{GShen}, Goncharov and Shen provide an elegant discussion of canonical bases from a geometric perspective.
Their discussion is similar to the approach developed in \cite{GHKK} (and employed here).
I'd like to point out key differences.
The results of Goncharov and Shen rely ultimately on
Lusztig. 
They do not produce a basis, they produce a set which they show (by
a very elegant argument) is in natural bijection with the
basis Lusztig constructed.
The main point of \cite{GHKK} is the construction of the basis.
From the construction, the fact that this basis is
parametrized by a canonical subset of the tropical points of the  
mirror is immediate from the start. 
Each choice of seed identifies this set with the  
integer points of a polyhedral cone.
A result established here is that for $U$ and $G/U$ there is a seed
in which the cone produced is already known and loved-- 
but the cones described would have corresponded to canonical bases independently of this fact.
Another difference is that the construction in \cite{GShen} is tailored to representation theory,
and it utilizes sophisticated representation theoretic tools like 
perverse sheaves on the affine Grassmannian and the geometric Satake correspondence.
On the other hand, the machinery of \cite{GHKK} is expected to apply to all affine log Calabi-Yau varieties with maximal boundary.
The results here are obtained essentially without representation
theoretic considerations. 
The fact that $\SL_n$ is a group, for example, is irrelevant to these constructions.
All that (conjecturally) matters is, for instance, that $G^{e,w_0}\subset G/U$ is an affine  
variety with a volume form of the right sort, which has a pole on each divisor in the
complement of $G^{e,w_0} \subset G/U$.
Many spaces of interest in representation theory have the requisite structure,
and this allows some results of \cite{GShen} to be obtained in a simpler fashion.}}

{\it{Acknowledgments: }}
My advisor, Sean Keel, introduced me to the topics discussed in this paper and suggested the questions I address here.
In addition to his guidance, I benefitted greatly from discussions with A. Fenyes, B. Leclerc, M. Lingam, T. Mandel, L. Shen, S. Villar, and H. Williams.
I also enjoyed RTG funding while writing this paper.

\section{Maximal unipotent in $\SL_n$}
The goal of this section is to establish the full Fock-Goncharov conjecture\footnote{%
This terminology is due to \cite{GHKK}, and they give a precise statement of the full Fock-Goncharov 
conjecture in Theorem 0.2 and Definition 0.6.
}
for the unipotent radical $U$ of the Borel subgroup $B_+$ of upper triangular matrices in $\SL_n$.\footnote{%
$U$ is just the subgroup of uppertriangular matrices with 1's along the diagonal.}
The space $U$ is of representation theory interest for the following reasons:
\begin{enumerate}
    \item Let $\mf{u}$ be the Lie algebra of $U$, and $\mathcal{U}_q(\mf{u})$ its quantized universal enveloping algebra.
	  Lusztig constructed his celebrated canonical basis in \cite{Lusztig_canonical_bases} precisely for $\mathcal{U}_q(\mf{u})$.\footnote{%
		In more detail, the paper deals with the positive part of a quantized enveloping algebra arising from a simply laced root system.  $\SL_n$ is of Cartan type $A_{n-1}$, and the positive part of the quantized enveloping algebra of its Lie algebra is $\mathcal{U}_q(\mf{u})$.%
}
	  It will be interesting to compare the canonical basis obtained via the machinery of cluster varieties and mirror symmetry to 
the classical limit of Lusztig's canonical basis.
    \item Now let $\mathcal{U}\lrp{\mf{u}}^*_{\gr}$ be the graded dual of the universal enveloping algebra of $\mf{u}$.
	  $H^0\lrp{U,\ssO_U}$ is isomorphic (as a Hopf algebra) to $\mathcal{U}\lrp{\mf{u}}^*_{\gr}$.\cite{GLS08}
    \item The base affine space $\SL_n/U$ has a very similar cluster structure to $U$. 
	  Analysis of $U$ in this section will provide the groundwork for analyzing $\SL_n/U$ 
	  in the next section, and some results will be proved simultaneously for both spaces.  
	  The coordinate ring of $\SL_n/U$ 
	  decomposes as a direct sum of (isomorphism classes of) irreducible representations of $\SL_n$.  
	  The canonical basis for this space consists of eigenfunctions of the natural torus action, 
	  so we will get a basis for each weight space by restriction.
	  That is, we will obtain a canonical basis on each irreducible representation of $\SL_n$.
\end{enumerate}

The first task is to establish $U$'s cluster structure.
We will see that $U$ arises as the partial compactification of an $\cA$-cluster variety, where all frozen variables are allowed to go to 0.
Some background sections are provided for completeness.
Clear and concise discussions of cluster varieties via initial data
are provided in Section 1.2 of \cite{FG_cluster_ensembles}
and Section 2 of \cite{GHK_birational}.
I won't reproduce the discussion here, but I will assume familiarity with it.
My notation follows \cite{GHK_birational}.
Feel free to skim the background sections for the results.

\subsection{Overview of double Bruhat cells, double pseudoline arrangements, and quivers for $G=\SL_n$}
This overview follows \cite{FZ_double_Bruhat} and \cite{FZ_clustersIII}. 

For $u,v$ in the Weyl group $W$ of $G$,
and $B_+, B_-$ a pair of opposite Borel subgroups,
the double Bruhat cell $G^{u,v}$ is defined to be
\eqn{
G^{u,v} = \lrp{B_+ u B_+} \cap \lrp{B_- v B_-}. 
}
Let $G = \SL_{n}$, $B_+$ the subgroup of upper triangular matrices, and $B_-$ lower triangular matrices.
$G$ is of type $A_{n-1}$ and $W = S_{n}$.

So the double Bruhat cells are indexed by $W\times W$.
Define a right action of the first copy of $W$ on $\lrc{1_R, 2_R, \dots, n_R}$,\footnote{%
This will be an indexing set for $n$ red pseudolines.
}
and a left action of the second copy of $W$ on $\lrc{1_B, 2_B, \dots, n_B}$.\footnote{%
This will be an indexing set for $n$ blue pseudolines.
}
\eqn{ 
j_R \cdot u = u\lrp{j_R} \qquad \qquad v \cdot j_B = v\lrp{j_B}
}
A {\it{double reduced word}} 
${\bf{i}}$
for $\lrp{u,v}$ is a shuffle of a reduced word   
${\bf{i}}_R$ for $u$ in the alphabet 
$\lrc{1_R, 2_R, \dots, (n-1)_R}$
and a reduced word 
${\bf{i}}_B$ for $v$ in the alphabet 
$\lrc{1_B, 2_B, \dots, (n-1)_B}$.

Representing ${\bf{i}}_R$ as a sequence of crossings of $n$ red ``pseudolines'' indexed $\lrc{1_R, 2_R, \dots, n_R}$ by height bottom up, 
$j_R$ corresponds to a crossing on the $j$-th level-- a crossing of the pseodolines in positions $j$ and $j+1$.
Similarly, $j_B$ represents a crossing of blue pseudolines on level $j$.
As $u$ acts on the right, the red pseudolines are indexed $\lrc{1_R, 2_R, \dots, n_R}$ on the left.
The blue pseudolines are labeled on the right.
These diagrams are called {\it{pseudoline arrangements}}.
As ${\bf{i}}$ is a shuffle of ${\bf{i}}_R$ and ${\bf{i}}_B$,
it can be represented by overlaying the pseudoline arrangements corresponding to ${\bf{i}}_R$ and ${\bf{i}}_B$.
The result is called a {\it{double pseudoline arrangement}}. 
For each ${\bf{i}}$,
there is a seed whose variables correspond to the chambers of the double pseudoline arrangement associated to ${\bf{i}}$.
For example, ${\bf{i}} = \lrp{ 1_B,1_R, 2_B, 1_B }$ produces the double pseudoline arrangement and chamber minors shown below.
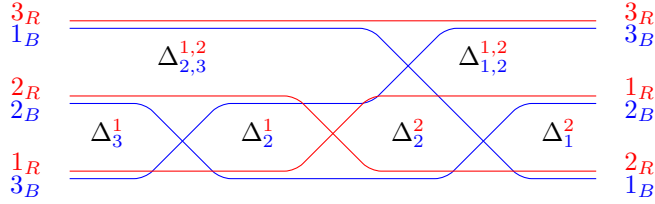
\begin{center}
\begin{tikzpicture}
	\draw[rounded corners, color = blue] (0,1.9) -- (1,1.9) -- (2,1.9) -- (3,1.9) -- (4,1.9) -- (5, 0.9) -- (6,-0.1) -- (7,-0.1) ;
	\draw[rounded corners, color = blue] (0,0.9) -- (1,0.9) -- (2,-0.1) -- (3,-0.1) -- (4,-0.1) -- (5, -0.1) -- (6,0.9) -- (7,0.9) ;
	\draw[rounded corners, color = blue] (0,-0.1) -- (1,-0.1) -- (2,0.9) -- (3,0.9) -- (4,0.9) -- (5, 1.9) -- (6,1.9) -- (7,1.9) ;
	\draw[rounded corners, color = red] (0,2) -- (1,2) -- (2,2) -- (3,2) -- (4,2) -- (5,2) -- (6,2) -- (7,2) ;
	\draw[rounded corners, color = red] (0,1) -- (1,1) -- (2,1) -- (3,1) -- (4,0) -- (5,0) -- (6,0) -- (7,0) ;
	\draw[rounded corners, color = red] (0,0) -- (1,0) -- (2,0) -- (3,0) -- (4,1) -- (5,1) -- (6,1) -- (7,1) ;

	\node[color = red, anchor=east] at (-0.25,2.1) {$3_R$};
	\node[color = red, anchor=west] at (7.25,2.1) {$3_R$};
	\node[color = red, anchor=east] at (-0.25,1.1) {$2_R$};
	\node[color = red, anchor=west] at (7.25,1.1) {$1_R$};
	\node[color = red, anchor=east] at (-0.25,0.1) {$1_R$};
	\node[color = red, anchor=west] at (7.25,0.1) {$2_R$};

	\node[color = blue, anchor=east] at (-0.25,1.8) {$1_B$};
	\node[color = blue, anchor=west] at (7.25,1.8) {$3_B$};
	\node[color = blue, anchor=east] at (-0.25,0.8) {$2_B$};
	\node[color = blue, anchor=west] at (7.25,0.8) {$2_B$};
	\node[color = blue, anchor=east] at (-0.25,-0.2) {$3_B$};
	\node[color = blue, anchor=west] at (7.25,-0.2) {$1_B$};

	\node at (0.5,0.5) {$\Delta^{\textcolor{red}{1}}_{\textcolor{blue}{3}}$};
	\node at (1.5,1.5) {$\Delta^{\textcolor{red}{1,2}}_{\textcolor{blue}{2,3}}$};
	\node at (2.5,0.5) {$\Delta^{\textcolor{red}{1}}_{\textcolor{blue}{2}}$};
	\node at (4.5,0.5) {$\Delta^{\textcolor{red}{2}}_{\textcolor{blue}{2}}$};
	\node at (5.5,1.5) {$\Delta^{\textcolor{red}{1,2}}_{\textcolor{blue}{1,2}}$};
	\node at (6.5,0.5) {$\Delta^{\textcolor{red}{2}}_{\textcolor{blue}{1}}$};

\end{tikzpicture}
\end{center}

The red
lines below a chamber index the rows for the minor $\Delta^{\textcolor{red}{I}}_{\textcolor{blue}{J}}$,
while the blue
lines below it index the columns.
All such seeds in $G^{u,v}$ are mutation equivalent.\cite[Remark 2.14]{FZ_clustersIII}

Berenstein, Fomin, and Zelevinsky give an explicit definition of the quiver associated to this seed, where chambers are the vertices of the quiver.\cite[Definition 2.2]{FZ_clustersIII}
I paraphrase their procedure below.
It is easiest to state without actively avoiding arrows between frozen variables.
Since these arrows do not affect exchange relations, I won't bother leaving them out.

Prepend the sequence $\lrp{(n-1)_R, (n-2)_R, \dots, 1_R} $ to ${\bf{i}}$ and call the result $\tilde{\bf{i}}$.
This introduces a sequence of fictitious red crossings on the left side of the diagram.

At a blue crossing, draw an arrow from the chamber at the left of the crossing to the chamber at the right.
Reverse the direction for a real red crossing.\footnote{%
This is a restatement of Definition 2.2, Condition (1) of \cite{FZ_clustersIII}.}
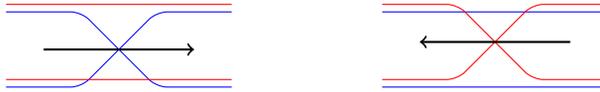
\begin{center}
\begin{tikzpicture}
	\draw[rounded corners, color = blue] (0,0.9) -- (1,0.9) -- (2,-0.1) -- (3,-0.1) ;
	\draw[rounded corners, color = blue] (0,-0.1) -- (1,-0.1) -- (2,0.9) -- (3,0.9) ;
	\draw[rounded corners, color = red] (0,1) -- (1,1) -- (2,1) -- (3,1) ;
	\draw[rounded corners, color = red] (0,0) -- (1,0) -- (2,0) -- (3,0) ;

	\draw[->,thick] (0.5,0.4) -- (2.5,0.4) ;

	\draw[rounded corners, color = blue] (5,0.9) -- (6,0.9) -- (7,0.9) -- (8,0.9) ;
	\draw[rounded corners, color = blue] (5,-0.1) -- (6,-0.1) -- (7,-0.1) -- (8,-0.1) ;
	\draw[rounded corners, color = red] (5,1) -- (6,1) -- (7,0) -- (8,0) ;
	\draw[rounded corners, color = red] (5,0) -- (6,0) -- (7,1) -- (8,1) ;

	\draw[<-,thick] (5.5,0.5) -- (7.5,0.5) ;
\end{tikzpicture}
\end{center}

Consider levels $j$ and $j'=j \pm 1$.\footnote{%
This is the condition $A_{\lrm{i_k},\lrm{i_l}} <0$ in \cite{FZ_clustersIII}.
}
Suppose there is exactly 1 crossing on level $j'$ between 2 consecutive crossings on level $j$.
Suppose further that some pseudoline is used in both the crossing on level $j'$ and the first crossing on level $j$, 
{\it{i.e.}} the crossing on level $j'$ 
and the first crossing on level $j$ are the same color.
If this color is red, draw an inclined arrow across this pseudoline pointing rightward.
If the pseudoline in question is blue, the inclined arrow will point leftward.\footnote{%
This is Definition 2.2, Condition (2) in \cite{FZ_clustersIII}.
}
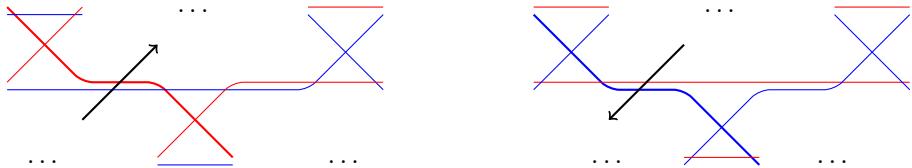
\begin{center}
\begin{tikzpicture}
	\draw[rounded corners, color = red, thick] (0,2) -- (1,1) -- (2,1) -- (3,0) ;
	\node at (2.5, 1.95) {$\dots$};
	\draw[rounded corners, color = blue] (0,1.9) -- (1,1.9) ;
	\draw[rounded corners, color = blue] (4,1.9) -- (5,0.9) ;
	\draw[rounded corners, color = red] (0,1) -- (1,2) ;
	\draw[rounded corners, color = red] (4,2) -- (5,2) ;
	\draw[rounded corners, color = blue] (0,0.9) -- (4,0.9) -- (5,1.9) ;
	\draw[rounded corners, color = blue] (2,-0.1) -- (3,-0.1) ;
	\draw[rounded corners, color = red] (2,0) -- (3,1) -- (5,1) ;
	\node at (0.5,-0.05) {$\dots$};
	\node at (4.5,-0.05) {$\dots$};

	\draw[->,thick] (1,0.5) -- (2,1.5) ;

	\draw[rounded corners, color = blue, thick] (7,1.9) -- (8,0.9) -- (9,0.9) -- (10,-0.1) ;
	\draw[rounded corners, color = red] (7,2) -- (8,2) ;
	\draw[rounded corners, color = red] (11,2) -- (12,2) ;
	\draw[rounded corners, color = red] (7,1) -- (8,1) -- (9,1) -- (10,1) -- (11,1) -- (12,1) ;
	\draw[rounded corners, color = blue] (7,0.9) -- (8,1.9) ;
	\draw[rounded corners, color = blue] (11,1.9) -- (12,0.9) ;
	\node at (9.5,1.95) {$\dots$};
	\draw[rounded corners, color = blue] (9,-0.1) -- (10,0.9) -- (11,0.9) -- (12,1.9) ;
	\draw[rounded corners, color = red] (9,0) -- (10,0) ;
	\node at (8,-0.05) {$\dots$};
	\node at (11,-0.05) {$\dots$};

	\draw[->,thick] (9,1.5) -- (8,0.5) ;
\end{tikzpicture}
\end{center}
Now suppose that following a crossing on level $j$, there is a pair of consecutive crossings on level $j'$ having opposite colors 
prior to any more crossings on $j$.
If the first of this pair of crossings is blue,
draw an arrow from the chamber bounded by the pair to the chamber initiated by the crossing on $j$.
If the colors for the pair of crossings are reversed, so is the direction of the arrow.\footnote{%
This is Definition 2.2, Condition (3)  of \cite{FZ_clustersIII}.
}
\begin{center}
\begin{tikzpicture}
	\draw[rounded corners, color = blue] (0,1.9) -- (1,0.9) -- (2,0.9) -- (3,-0.1) -- (4,-0.1) -- (5,-0.1) ;
	\draw[rounded corners, color = red] (0,2) -- (1,2) ;
	\draw[rounded corners, color = red] (0,1) -- (1,1) -- (2,1) -- (3,1) -- (4,1) -- (5,0);
	\draw[rounded corners, color = blue] (0,0.9) -- (1,1.9) ;
	\node at (2.5,1.95) {$\dots$};
	\draw[rounded corners, color = blue] (2,-0.1) -- (3,0.9) -- (4,0.9) -- (5,0.9) ;
	\draw[rounded corners, color = red] (2,0) -- (4,0) -- (5,1) ;
	\node at (0.5,-0.05) {$\dots$};

	\draw[->,thick] (3.5,0.5) -- (2.5,1.5) ;

	\node at (7.5,1.95) {$\dots$};
	\draw[rounded corners, color = blue] (9,1.9) -- (11,1.9) -- (12,0.9) ;
	\draw[rounded corners, color = red] (9,2) -- (10,1) -- (12,1) ;
	\draw[rounded corners, color = red] (7,1) -- (9,1) -- (10,2) -- (12,2) ;
	\draw[rounded corners, color = blue] (7,0.9) -- (8,-0.1) ;
	\node at (9.5,-0.05) {$\dots$};
	\draw[rounded corners, color = blue] (7,-0.1) -- (8,0.9) -- (11,0.9) -- (12,1.9) ;
	\draw[rounded corners, color = red] (7,0) -- (8,0) ;

	\draw[->,thick] (9.5,0.5) -- (10.5,1.5) ;
\end{tikzpicture}
\end{center}
Note that arrows only connect neighboring chambers,
and these arrows are fully determined by the pseudolines bounding the pair of chambers in question.
 



\subsection{$U$ as partial compactification of an $\cA$-cluster variety}

The unipotent radical $U$ of $B_+$ is the subgroup of upper triangular matrices with 1's along the diagonal.
Consider the double Bruhat cell $G^{e,w_0} = \lrp{B_+ e B_+}\cap \lrp{B_- w_0 B_-}$.
The intersection $\mathring{U} := U \cap G^{e,w_0}$ is the non-vanishing locus of the upper right minors in $U$.\footnote{See Fomin and Zelevinsky's description of Bruhat cells in terms of generalized minors in \cite{FZ_double_Bruhat}.} 
$\mathring{U}$ is an $\cA$-cluster variety, and $U$ is obtained from $\mathring{U}$ by allowing all frozen variables on $\mathring{U}$ to vanish.

The usual initial data for $\mathring{U}$ comes from initial data for $G^{e,w_0}$.
The initial data generally chosen for $G^{e,w_0}$ is associated to the reduced word for $w_0$ that is minimal in the lexicographic ordering.
For example, if $n=5$, we would choose the expression
$\lrp{1,2,1,3,2,1,4,3,2,1}$
for $w_0$.
The resulting double pseudoline arrangement and chamber minors are shown below.
\input{DPLA.tex}

The minors associated to unbounded chambers are the frozen variables.
Note that each chamber minor on the right side of the double pseudoline arrangement is identically 1 when restricted to $\mathring{U} \subset G^{e,w_0}$.
We can simply ignore these variables and the associated vertices and arrows in the quiver if we are describing $\mathring{U}$ as a cluster variety.
The resulting quiver $Q$ for the $n=5$ example is shown below, with frozen vertices boxed.

\input{U_quiver.tex}

I'll write $A_{{i;j}_s}$ to denote the cluster variable associated to the vertex $v_{i;j}$ in the seed $s$, or just $A_{i;j}$ if the seed is clear from context.
I'll write $I$ for the indexing set of the vertices of the quiver, or equivalently the indexing set of the variables in a given seed.

So choosing the lexicographically minimal reduced word $\vb{i}_{w_0}$ for $w_0$ determines an open torus in $\mathring{U}$ 
as the non-vanishing locus of the chamber minors associated to $\vb{i}_{w_0}$.
This is the torus associated to the initial seed $s_0$ in the toric atlas that defines $\mathring{U}$'s cluster variety structure.
However, there is another natural way in which $\vb{i}_{w_0}$ determines a torus in $\mathring{U}$.
For a given semisimple Lie group $G$, we can define the following one parameter subgroups of $G$ via the Chevalley-Serre generators for the Lie algebra.
\eqn{
E_i(t):= \exp\lrp{t e_i}\\
F_i(t):= \exp\lrp{t f_i}\\
H_i(t):= \exp\lrp{t h_i}
}

Note that only the $E_i(t)$ subgroups land in $U$.
The choice of reduced word $\vb{i}_{w_0}$ determines the following immersion initially described by A. Whitney and Loewner in the context of totally positive matrices, then generalized to the setting of reductive groups and related to canonical bases by Lusztig:\cite{AWhitney,Loewner,Lusztig_TPRed,Lusztig_TPCB}
\eqn{
\varphi_{\vb{i}_{w_0}}:\lrp{\C^*}^{l\lrp{w_0}} &\to U\\
\lrp{t_1, \dots, t_{l\lrp{w_0}}} &\mapsto E_{i_1}\lrp{t_1} \cdots E_{i_{l\lrp{w_0}}}\lrp{t_{l\lrp{w_0}}},
}
where $l\lrp{w_0}$ is the length of a reduced expression for $w_0$.

So $\vb{i}_{w_0}$ has determined two tori inside of $U$-- 
one as the non-vanishing locus of chamber minors 
and the other via $\varphi_{\vb{i}_{w_0}}$.
It is natural to wonder if these tori are related.

{\prop{\label{prop:reflect} Let $\alpha: U \to U$ be the automorphism of $U$ given by reflection over the anti-diagonal.  
Then the image of $\varphi_{\vb{i}_{w_0}}$ is the non-vanishing locus of $\lrc{\alpha^*\lrp{A_{{i;j}_{s_0}}}}_{i;j \in I}$, {\it{i.e.}} the non-vanishing locus of $\lrc{\Delta^{w_0 (J)}_{w_0(I)}}$ where each $\Delta^I_J$ is an initial seed minor.}}

\begin{proof}
Using Proposition 4.2 of \cite{FZ_double_Bruhat}, the following {\it{planar network}} can be used to write down $\varphi_{\vb{i}_{w_0}} \lrp{t_1, \dots, t_{l\lrp{w_0}} }$.

\input{phi_diagram.tex}

All paths are direct from left to right.
The weight of a path is the product of the weights of its segments,
and 
$\Delta^i_j \lrp{\varphi_{\vb{i}_{w_0}} \lrp{t_1, \dots, t_{l\lrp{w_0}}}}$ 
is the sum of the weights of all paths from $i$ to $j$. 
For example, for $n=4$,\footnote{In this case the ellipses in the picture can be ignored.} 
\eqn{
\Delta^1_3\lrp{\varphi_{\vb{i}_{w_0}} \lrp{t_1, \dots, t_6}} = t_1 t_2 + t_1 t_5 + t_3 t_5,
}
\eqn{
\Delta^3_3\lrp{\varphi_{\vb{i}_{w_0}} \lrp{t_1, \dots, t_6}} = 1,
}
and
\eqn{
\Delta^4_3\lrp{\varphi_{\vb{i}_{w_0}} \lrp{t_1, \dots, t_6}} = 0.
}
Now consider the following relabeling of the $t_i$'s:
\eqn{
\tau_{j;k} = t_{l} \iff k-1 = w_0\lrp{i_l} \text{ and the $l^{\text{th}}$ term in $\vb{i}_{w_0}$ is the $j^{\text{th}}$ occurrence of $i_l$ in $\vb{i}_{w_0}$}.
}
For example, for $n=4$ Figure {\ref{fig:planar}} becomes 
 
\input{replanar.tex}

Now define 
\eqn{
a_{i;j} = \begin{cases} 
\alpha^*\lrp{A_{{i;j}_{s_0}}} &\text{if $i;j$ is in $I$}\\
1 &\text{otherwise.}
\end{cases}
}
I claim that 
\eqn{
a_{i;j} \lrp{ \varphi_{\vb{i}_{w_0}} \lrp{ t_1, \dots, t_{l\lrp{w_0} }}} = \prod_{k=1}^i \prod_{l=k+1}^{k+j-i} \tau_{k;l}.
}
Note that 
$\lrc{ \prod_{k=1}^i \prod_{l=k+1}^{k+j-i} \tau_{k;l}}_{i;j \in I}$
are all non-zero precisely when
$\lrc{\tau_{i;j}}_{i;j \in I}$ are all non-zero.
Furthermore, if I order the $a_{i;j}$ lexicographically,
the product for $a_{i;j}$ has exactly 1 factor that has not appeared in a previous product-- namely, $\tau_{i;j}$.
So, assuming the claim, it is rather easy to explicitly construct an inverse 
$$\psi:  \lrp{\C^*}^{l\lrp{w_0}}_{{a_{1;2}}^{\pm 1}, {a_{1;3}}^{\pm 1}, \dots, {a_{n-1;n}}^{\pm 1}} \to \lrp{\C^*}^{l\lrp{w_0}}_{{\tau_{1;2}}^{\pm 1},{\tau_{1;3}}^{\pm 1}, \dots, {\tau_{n-1;n}}^{\pm 1} } $$ 
to $\varphi_{\vb{i}_{w_0}}$.
Then the proposition follows immediately from this claim.

Essentially, my proof of the claim will make use of multilinearity of the determinant to ``factor out a weight'' of a segment
and skew symmetry to ``delete segments''
until all that remains is a product consisting of the desired weights and the determinant of an upper triangular matrix with 1's on the diagonal.
Since
\eqn{
a_{i;j} = \Delta^{w_0(j), \dots, w_0(j -i +1)}_{w_0(i), \dots, w_0(1)},
}
the relevant submatrix $M$ is the $i\times i$ right aligned block
whose bottom row is $w_0\lrp{j-i+1}$.
So the paths relevant to the last column start between $w_0(j)$ and $w_0(j-i+1)$
and end at $w_0(1)$,
the paths relevant to the penultimate column start between $w_0(j)$ and $w_0(j-i+1)$
and end at $w_0(2)$,
and so forth.
This is illustrated below for $a_{2;4}$.

\input{a24paths.tex}

Note that every path used to write the last column of $M$ contains the line segment with weight $\tau_{1;2}$.
This indicates that $\tau_{1;2}$ is divides every entry of this column.
By multilinearity, this can be factored out of the last column yielding a new matrix $M'$ with $\det (M) = \tau_{1;2} \det (M')$.
The paths describing the entries of last column of $M'$ have the same weights as paths 
describing the entries in the same row of the second to last column.\footnote{%
Assuming there is more than 1 column-- and the result can be read off of the planar network immediately if there is not more than 1 column.
}
For example, the paths indicated below have the same weight.

\input{same_weight.tex}
  
So the last column of $M'$ is a summand of the second to last column.
By skew symmetry, this summand can be removed from the second to last column without affecting $\det (M')$.
This corresponds to deleting a segment of the planar network, as shown below.

\input{del_segment.tex}

As a result of deleting this segment, each path now used to describe the second to last column contains the segment
of weight $\tau_{2;3}$.
So this can be factored out the same way $\tau_{1;2}$ was previously.
If there are more than 2 columns, factoring out $\tau_{2;3}$ and using skew symmetry of the determinant allows me to 
delete the horizontal segment immediately below $\tau_{2;3}$'s segment.
(If there aren't more than 2 columns, the segment in question is never used.  Deleting it won't harm things; it's just unnecessary.)
Meanwhile, if $j-i>1$, then every path describing the last column contains the segment of weight $\tau_{1;3}$.
Again, factor it out and delete the horizontal segment below it, now killing paths used to describe both the second and third to last column.
This process continues until the paths for each column no longer share a segment with weight different from 1,
{\it{i.e.}} until the bottom non-zero entry of each column in the $i\times i$ matrix is 1.
This occurs after $j-i$ factors have been pulled off of the paths for each of the $i$ columns-- the factors for the last column having first index 1, those for the second to last having first index 2, and so forth--  
leaving a single weight 1 path from $w_0(j-i+k)$ to $w_0(k)$.
Continuing with the $a_{2;4}$ example, the resulting planar network has the form

\input{a24done.tex}

In this example, the weights 
$\tau_{1;2}$, $\tau_{1;3}$, $\tau_{2;3}$, and $\tau_{2;4}$ have been factored out, 
and the submatrix described by the final planar network is 
$\matdd{1}{\tau_{1;4}}{0}{1}$.
More generally, we have arrived at
\eqn{
\det(M) = \det(U) \prod_{k=1}^i \prod_{l = k+1}^{j-i+k} \tau_{k;l} ,
}
where $U$ is an upper triangular matrix with 1's on the diagonal.
The claim follows.
\end{proof}

In fact, the relationship between these two tori goes deeper.
{\prop{\label{prop:2torimutequiv}
The image of $\varphi_{\vb{i}_{w_0}}$ is also in the toric atlas for $\mathring{U}$.  The seed $s'$ for this torus can be reached from the initial seed $s_0$ by the following sequence of mutations:

$$v_{1;2}, v_{1;3}, v_{2;3}, \dots, v_{n-2;n-1},v_{1;2}, v_{1;3}, v_{2;3}, \dots, v_{n-3;n-2}, \dots , v_{1;2},v_{1;3},v_{2;3},v_{1;2} . $$

Furthermore, if $\mu$ is the resulting birational map $\mu: T_{N;s_0} \dashrightarrow T_{N;s'}$, 
then $\mu^*\lrp{A_{{i;j}_{s'}}} = \alpha^* \lrp{A_{{i;j}_{s_0}}}$. 
The quiver $Q_{s'}$ for $s'$ is just the quiver $Q_{s_0}$ for $s_0$ with every arrow reversed.
}}

The sequence listed in the proposition is easier to digest with a picture, 
so before I prove the proposition, I'll illustrate this sequence on the $n=5$ quiver from Figure~\ref{fig:U_quiver}.
Replacing vertex labels with their order in the sequence, the first 6 mutations are:

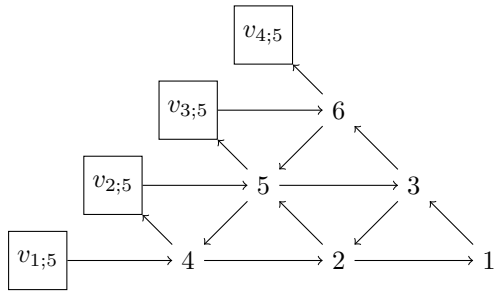
\input{vert_seq1.tex}

This is followed by

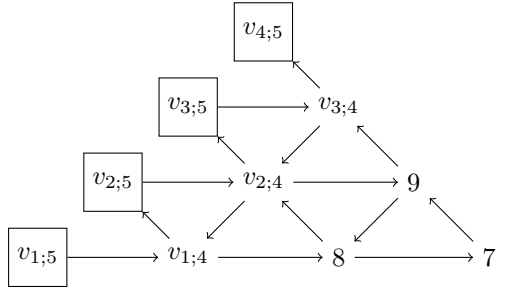
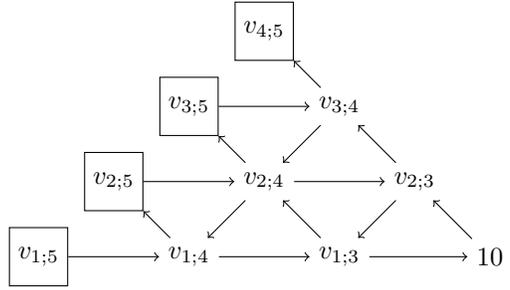
\input{vert_seq2.tex}

and finally

\input{vert_seq3.tex}

I'll prove Proposition~\ref{prop:2torimutequiv} in several steps.
Note that by Proposition~\ref{prop:reflect}, if $\mu^*\lrp{A_{{i;j}_{s'}}} = \alpha^* \lrp{A_{{i;j}_{s_0}}}$,
then $T_{N;s'} = \varphi_{\vb{i}_{w_0}}\lrp{\lrp{\C^*}^{l\lrp{w_0}}} $.
So to prove Proposition~\ref{prop:2torimutequiv}, 
it is sufficient to show that the listed sequence of mutations produces the claimed variables and quiver.
The quivers for each seed will be needed to mutate the variables, so I'll prove the claim about quivers first.
The following lemma will be useful.
{\lemma{\label{lem:delbot}
Let $Q_n$ be the quiver

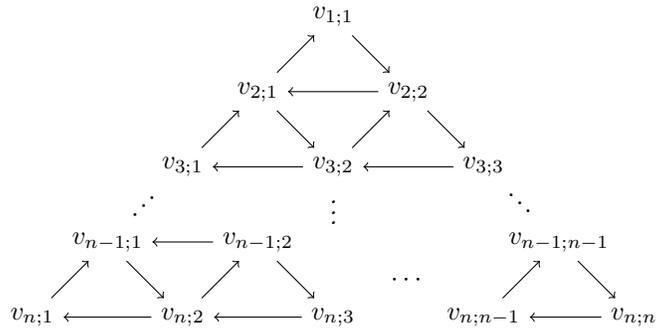
\input{Qn.tex}

Mutating at each vertex of $Q_n$ in lexicographic order simply deletes the bottom row of arrows, leaving

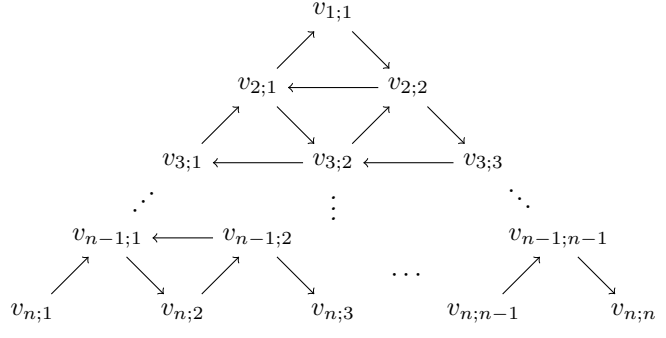
\input{Qnp.tex}

Denote this resulting quiver $\widehat{Q_n}$.
}}
\begin{proof}
For $Q_1$, there is nothing to show.
$Q_2$ mutates as follows:

\input{Q2mut.tex}

Now suppose it holds for all $m <n$.
Identify $Q_m$ with the subquiver of $Q_n$ containing the first $m$ rows.
Then the quiver produced by mutating $Q_{n-1}$ at every vertex in $Q_{n-2}$ is 
given by mutating $\widehat{Q_{n-1}}$ at $v_{n-1;n-1}, v_{n-1;n-2}, \dots, v_{n-1;1}$.
(This is just inverting the last row of mutations.)
The relevant portion of the quiver is

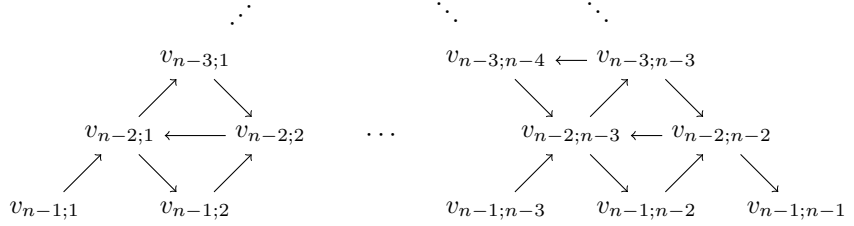
\input{Qnm1hatmut.tex}

Mutating at $v_{n-1;n-1}$ yields

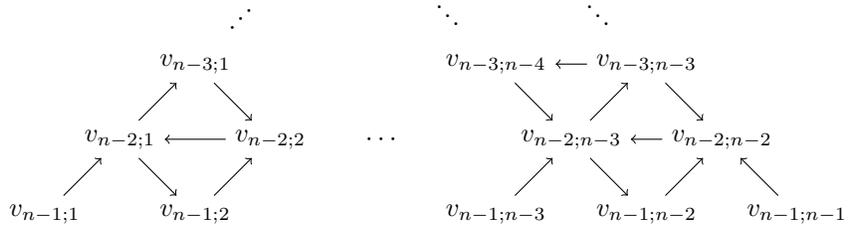
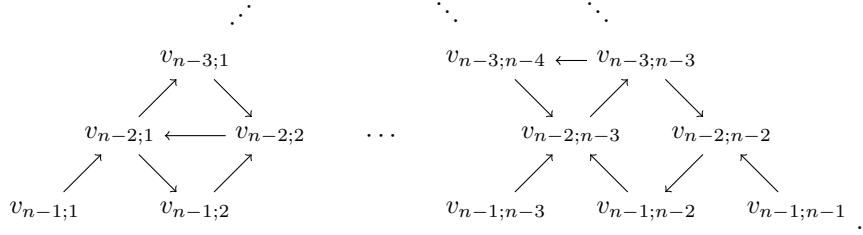
\input{Qnm1hatmut1.tex}

And $v_{n-1;n-2}$:

\input{Qnm1hatmut2.tex}

As the sequence continues, so does this process of reversing the direction of arrows to and from 
the bottom row while cancelling horizontal arrows on row $n-2$.
The end result is 

\input{Qnm1hatmutdone.tex}

Now note that there are no arrows in $Q_n$ between row $n$ and $Q_{n-2}$.
As a result, no set of mutations that is restricted to $Q_{n-2}$ can affect 
arrows involving the $n^{\text{th}}$ row of vertices,
nor can the presence of the $n^{\text{th}}$ row affect the mutation of $Q_{n-2}$ in $Q_{n-1}$.
So, mutating through all of the vertices of $Q_{n-2}$ must produce the following quiver, 
which contains an embedded copy of the quiver in Figure~\ref{fig:Qnm1muted}:

\input{Qnafter2.tex}

A similar line of reasoning shows that none of the remaining mutations (all of which are on rows $n-1$ and $n$) can affect $Q_{n-3}$, or any of the arrows with a vertex in $Q_{n-3}$.
Then I can mutate the following more manageable subquiver instead:

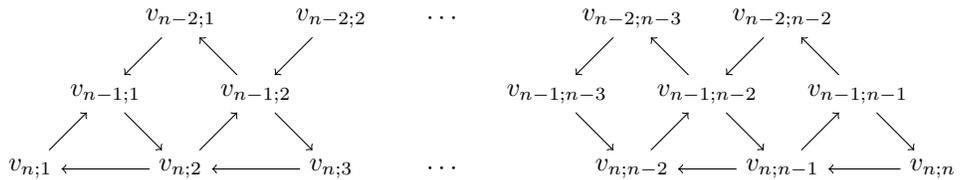
\input{Qnpbottom.tex}

Mutating at $v_{n-1;1}$ yields:

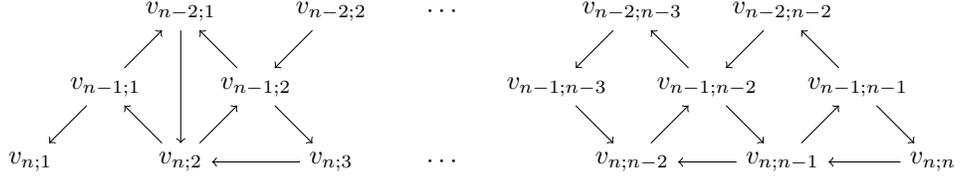
\input{Qnpbottom1.tex}

Following with a mutation at $v_{n-1;1}$ gives:

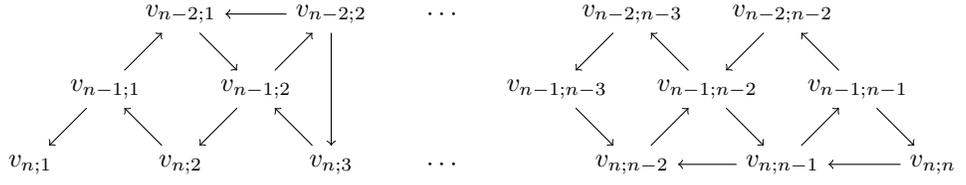
\input{Qnpbottom2.tex}

A clear pattern develops.
Mutating at $v_{n-1;k}$ kills a downward vertical arrow to its left while creating downward vertical arrow to its right.
Similarly, a leftward horizontal arrow below $v_{n-1;k}$ is replaced by a leftward horizontal arrow above it.
Meanwhile, the direction of each incoming or outgoing arrow at $v_{n-1;k}$ is reversed.
Then the result of mutating at each of the vertices in row $n-1$ is

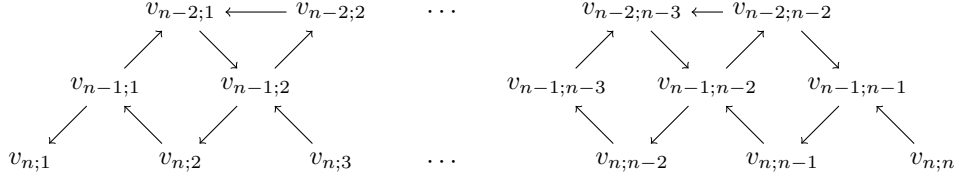
\input{Qnpbottomnm1.tex}

Now compare this quiver to Figure~\ref{fig:Qnm1muted}, 
and recall that Figure~\ref{fig:Qnm1muted} was constructed such that mutating along the bottom row produced $\widehat{Q_{n-1}}$.
This construction applied regardless of the number of rows in the quiver--
it did not make use of the induction hypothesis.
Instead, the induction hypothesis was used to say that Figure~\ref{fig:Qnm1muted} 
resulted from mutating $Q_{n-1}$ along its first $n-2$ rows.
Then mutating along the $n^{\text{th}}$ row produces the bottom of $\widehat{Q_n}$, proving the lemma.
\end{proof}

{\cor{\label{cor:revarrows}  
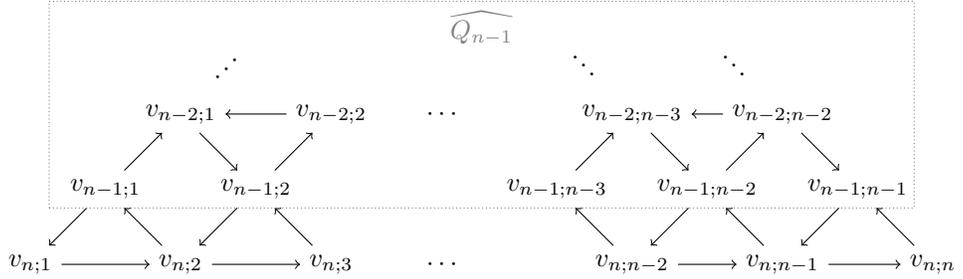
Mutating $Q_n$ at each vertex lexicographically, then mutating at each vertex of its subquiver $Q_{n-1}$, and so forth finishing with $Q_1$ simply reverses the direction of every arrow in $Q_n$.
}}

\begin{proof}
As before, for $Q_1$, there is nothing to show, and for $Q_2$
there is only one mutation to add to the sequence in Figure~\ref{fig:Q2mut}:

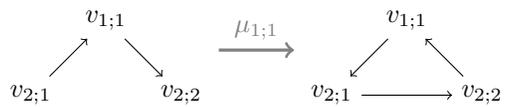
\input{Q2mut2.tex}

So suppose the claim holds for all $m <n$.
Mutating at each vertex of $Q_n$ yields $\widehat{Q_n}$ by Lemma~\ref{lem:delbot}.
Since there are no arrows between the subquiver $Q_{n-2}$ and the $n^\text{th}$ row of $\widehat{Q_n}$,
the argument used to draw Figure~\ref{fig:Qnafter2} applies here as well.
Mutating $\widehat{Q_n}$ through the first $n-2$ rows produces

\input{Qnhatafter2.tex}

The quiver mutations following Figure~\ref{fig:Qnafter2}
transfer to this setting almost exactly.
The one difference is that the current quiver has no 3-cycles at the bottom.
Since there are no leftward horizontal arrows along the bottom row to cancel,
mutation will create rightward horizontal arrows on the bottom row.
So in the case, the result of mutating at the vertices of row $n-1$ is 

\input{Qnhatnm1.tex}

All arrows associated with the bottom row of vertices 
are now in the desired position.
None of the remaining mutations can affect these arrows.
Only the subquiver $\widehat{Q_{n-1}}$ will change.
But by Lemma~\ref{lem:delbot}, this subquiver is the result of mutating
$Q_{n-1}$ at each of its vertices in lexicographic order.
Then the remaining mutations turn 
Figure~\ref{fig:Qnmutnnm1}
into $Q_n$ with each arrow reversed by the induction hypothesis.
\end{proof}

{\cor{
The sequence of mutations listed in Proposition~\ref{prop:2torimutequiv} simply reverses every arrow of $Q_{s_0}$. 
}}
\begin{proof}
$Q_{s_0}$ is $\widehat{Q_n}$ with vertices relabeled.
The frozen vertices of $Q_{s_0}$ correspond to the bottom row of $\widehat{Q_n}$.
Under the relabeling, the listed sequence says to mutate at the vertices of $Q_{n-1}$ lexicographically, then $Q_{n-2}$, and so forth.
So this is precisely the case considered in Corollary~\ref{cor:revarrows} after mutating through the vertices of $Q_n$ lexicographically.
The resulting quiver is $Q_n$ with every arrow reversed, and deleting arrows between frozen vertices yields $Q_{s_0}$ with each arrow reversed.
\end{proof}

To prove Proposition~\ref{prop:2torimutequiv}, all that remains is mutating the variables of the initial seed.
I claim that each mutation at $v_{i;j}$ in the sequence shifts the indexing sets for both the rows and the columns up by 1,
{\it{e.g.}} $\lrc{3,4,5}$ is sent to $\lrc{4,5,6}$.
Note that if this holds, the sequence of mutations would reflect each $A_{{i;j}_{s_0}}$ over the antidiagonal, finishing the proof 
of Proposition~\ref{prop:2torimutequiv}. 
The proof of this claim will be another exercise in quiver mutations, together with repeated application of a certain minor identity.

{\lemma{\label{lem:minorIdent}
Consider the space $\Mat_{n \times (n+1)}\lrp{\mathbb{K}}$, where $\mathbb{K}$ is a field of characteristic $\neq 2$.
Let $J = \lrc{1,\dots, n+1}$ be the indexing set for columns. 
Choose $j_1$, $j_2$, and $j_3$, with $1 \leq j_1 < j_2 < j_3 \leq n+1$.
Then 
\eqn{
\Delta^{1,\dots, n-1}_{J-\lrc{j_2,j_3}}
\Delta^{1,\dots, n}_{J- \lrc{j_1}}
-
\Delta^{1,\dots, n-1}_{J- \lrc{j_1,j_3}}
\Delta^{1,\dots, n}_{J- \lrc{j_2}}
+
\Delta^{1,\dots, n-1}_{J- \lrc{j_1,j_2}}
\Delta^{1,\dots, n}_{J- \lrc{j_3}}
= 0
}
on $\Mat_{n \times (n+1)}$.
}}

\begin{proof}\footnote{Much of this proof is due to Aaron Fenyes.}
The minors in question are determined by choosing (with order) $n+1$ vectors in an $n$-dimensional vector space $V$,
with each vector representing a column of a matrix in $\Mat_{n\times (n+1)}$.
The vectors that do not correspond to columns $j_1$, $j_2$, or $j_3$ occur in each minor.
There are $n-2$ of these vectors. Fix them.
Then 
$$ 
\Delta^{1,\dots, n-1}_{J-\lrc{j_2,j_3}}
\Delta^{1,\dots, n}_{J- \lrc{j_1}}
-
\Delta^{1,\dots, n-1}_{J- \lrc{j_1,j_3}}
\Delta^{1,\dots, n}_{J- \lrc{j_2}}
+
\Delta^{1,\dots, n-1}_{J- \lrc{j_1,j_2}}
\Delta^{1,\dots, n}_{J- \lrc{j_3}}
$$
becomes a function $f$ of the remaining 3 vectors.
If the fixed vectors are not linearly independent, the $n \times n$ minors will vanish.
Then $f$ is identically 0, and the claim holds.
So assume the fixed vectors are linearly independent. 
Write $v_1$, $v_2$, and $v_3$
for the vectors corresponding to $j_1$, $j_2$, and $j_3$,
and write $v_i = v_i^\perp + v_i^\parallel$ 
to decompose these vectors into the components perpendicular and parallel to the subspace $W$ spanned by the fixed vectors.
$v_i^\parallel$ does not affect $f$, as it does not affect the minors defining $f$.
We are left with $v_1^\perp$, $v_2^\perp$, and $v_3^\perp$ in $W^\perp$.
$f$ is clearly multilinear in these arguments.
It is also skew-symmetric, as described below. 
If any two of $\lrc{v_1^\perp,v_2^\perp,v_3^\perp}$ are swapped, the sign of the $n \times n$ minor that utilizes both of these vectors is flipped by skew-symmetry of the determinant.  
The $(n-1) \times (n-1)$ minor paired with it is independent of the swapped vectors, so it remains unaltered.
So swapping the two vectors flips the sign of this summand.
Meanwhile, if the remaining summands are $s$ and $s'$, swapping the two vectors swaps $s$ and $s'$ up to a sign. 
This is because it exchanges the sets of vectors that determine the minors in these summands, without respecting the ordering on these sets.
For instance, upon swapping $v_1^\perp$ and $v_2^\perp$,
$\Delta^{1,\dots, n-1}_{J-\lrc{j_2,j_3}}$
becomes a signed volume of the $\lrp{n-1}$-dimensional parallelepiped 
whose edges are projections of the fixed vectors and $v_2$.\footnote{%
The edges are obtained by projecting out the last coordinate of these vectors.  The parallelepiped lives in the hyperplane of $V$ where the last coordinate is 0.
}
$\Delta^{1,\dots, n-1}_{J- \lrc{j_1,j_3}}$
is also a signed volume of this parallelepiped.
The relative sign of these two volumes is determined by the parity of the permutation sending one ordered set of edges to the other.
There are $j_2-j_1-1$ positions between $j_1$ and $j_2$,
so it takes $j_2-j_1-1$ neighboring transpositions to send the first ordered set to the second, 
and the relative sign is $(-1)^{j_2-j_1-1}$.
Similarly,
$\Delta^{1,\dots, n}_{J- \lrc{j_1}}$
becomes 
$
(-1)^{j_2-j_1-1} 
\Delta^{1,\dots, n}_{J- \lrc{j_2}}.
$
So if $s = \Delta^{1,\dots, n-1}_{J-\lrc{j_2,j_3}}\Delta^{1,\dots, n}_{J- \lrc{j_1}}$,
swapping $v_1^\perp$ and $v_2^\perp$ sends $s$ to 
\eqn{
&(-1)^{j_2-j_1-1} 
(-1)^{j_2-j_1-1} 
\Delta^{1,\dots, n-1}_{J- \lrc{j_1,j_3}}
\Delta^{1,\dots, n}_{J- \lrc{j_2}}\\
=&
\Delta^{1,\dots, n-1}_{J- \lrc{j_1,j_3}}
\Delta^{1,\dots, n}_{J- \lrc{j_2}}\\
=&-s'.
}
The same argument shows that this swapping of $v_1^\perp$ and $v_2^\perp$ sends $s'$ to $-s$ as well.
Swapping $v_1^\perp$ and $v_2^\perp$ has introduced an overall minus sign.
This argument transfers rather directly to swapping $v_2^\perp$ and $v_3^\perp$.
The only difference that arises when swapping $v_1^\perp$ and $v_3^\perp$ is
that the column indexed by $j_2$ is missing from the minors 
$\Delta^{1,\dots, n-1}_{J-\lrc{j_2,j_3}}$
and
$\Delta^{1,\dots, n-1}_{J- \lrc{j_1,j_2}}$,
so the relative sign of the relevant volumes will be $(-1)^{j_3 - j_1 -2}$ instead of $(-1)^{j_3 - j_1 -1}$.
This factor of $-1$ accounts for the fact that the coefficients of both
$\Delta^{1,\dots, n-1}_{J-\lrc{j_2,j_3}}
\Delta^{1,\dots, n}_{J- \lrc{j_1}}
$
and 
$
\Delta^{1,\dots, n-1}_{J- \lrc{j_1,j_2}}
\Delta^{1,\dots, n}_{J- \lrc{j_3}}
$
are $+1$ in the expression for $f$, in contrast to the differing signs of the previous two cases.
Then swapping any two of $\lrc{v_1^\perp,v_2^\perp, v_3^\perp}$ send $f$ to $-f$.
$f$ is a skew-symmetric multilinear form on $V^3$, 
but $v_1^\perp$, $v_2^\perp$, and $v_3^\perp$ all lie in the two dimensional space $W^\perp$.
A skew-symmetric multilinear form on $V^n$ vanishes on $n$-tuples of linearly dependent vectors,\footnote{Recall that $\chara(\mathbb{K}) \neq 2$.} so $f$ is again identically 0.
The identity holds.   
\end{proof}

Using Lemma~\ref{lem:minorIdent} repeatedly yields the following proposition for mutations in $G^{e,w_0}$.
{\prop{\label{prop:Gew0minors}
Consider the sequence of mutations in Proposition~\ref{prop:2torimutequiv}.
After the $k^\text{th}$ mutation at $v_{i;j}$,
$A_{i;j}$ pulls back to $T_{N;s_0}$ as $\Delta^{1,\dots,k,k+1,\dots, k+i}_{1,\dots,k, k+j-i+1,\dots, k+ j}$.
}}
\begin{proof}
While it was convenient to drop the frozen vertices $v_{i;i}$ before, it's necessary to reintroduce them now.
The differences between mutations of the new quivers and those investigated previously are uncomplicated,
so rather than repeating the previous discussion with minor changes, I'll simply sketch the results.
Instead of $\widehat{Q_n}$ from Lemma~\ref{lem:delbot},
the initial seed quiver $\widetilde{Q_{s_0}}$ is now a relabeling of

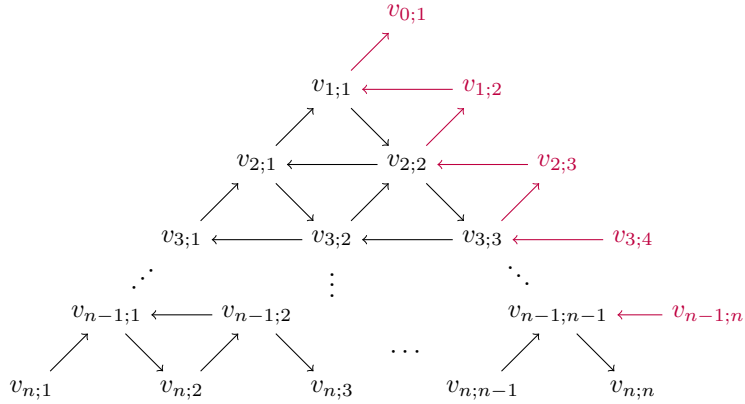
\input{Qtilden.tex}

The reintroduced parts of the quiver are colored for clarity. Anything colored did not appear in $\widehat{Q_n}$.
Call the new quiver $\widetilde{Q_n}$.
The bottom row and the new vertices are frozen, so I'll leave out arrows that would arise between these vertices.
The sequence of mutations only produces new arrows between the new vertices and vertices labeled $v_{i;1}$ or $v_{i;i}$.
Mutating through the vertices of $Q_{n-1}$ lexicographically produces

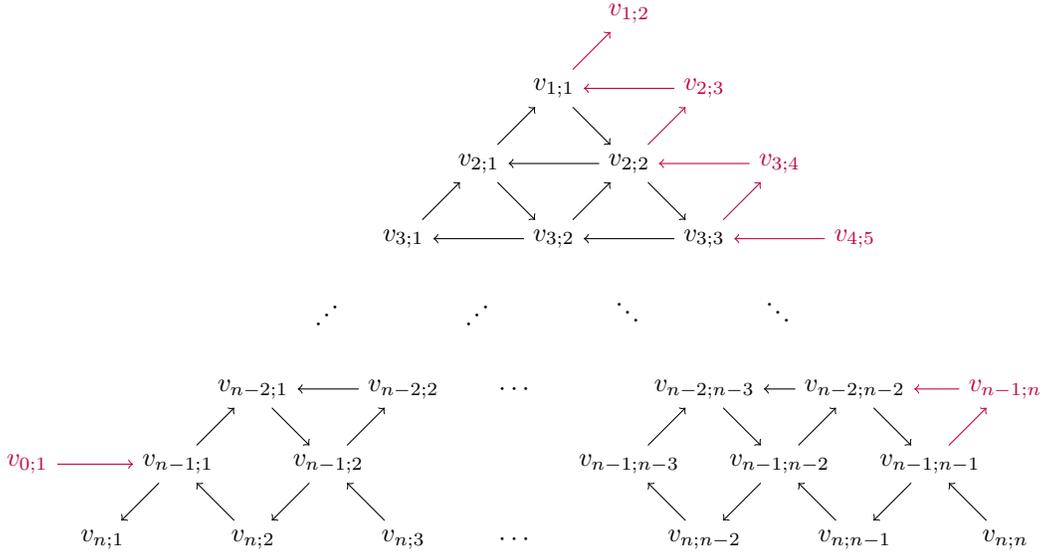
\input{Qtildenmut1.tex}

The full sequence of mutations yields

\input{Qtildenmutdone.tex}

Each mutation at $v_{i;j}$ utilizes a quiver in which $v_{i;j}$ has two outgoing arrows and two incoming arrows.
Specifically, there is an arrow from $v_{i+1;j}$ to $v_{i;j}$ and an arrow from $v_{i;j}$ to $v_{i+1;j+1}$.
If $v_{i-1;j-1}$ is a vertex of the subquiver $\widehat{Q_n}$ of $\widetilde{Q_n}$, then there is an arrow from $v_{i;j}$ to $v_{i-1;j-1}$.
If $v_{i-1;j}$ is a vertex of $\widehat{Q_n}$, then there is an arrow from $v_{i-1;j}$ to $v_{i;j}$.
If $v_{i-1;j-1}$ is not a vertex of $\widehat{Q_n}$, then $j=1$.
In this case, in the quiver for the $k^\text{th}$ mutation at $v_{i;j}$ there is an arrow from $v_{i;j}$ to $v_{k-1;k}$. 
If $v_{i-1;j}$ is not a vertex of $\widehat{Q_n}$, then $j=i$.
In the quiver for the $k^\text{th}$ mutation at $v_{i;j}$, there is an arrow from $v_{k+i-1;k+i}$ to $v_{i;j}$.

The relabeling of $\widetilde{Q_n}$ needed to obtain $\widetilde{Q_{s_0}}$ is $v_{i;j} \mapsto v_{j;i+1}$.
Then, using the labeling of $\widetilde{Q_{s_0}}$, 
the exchange relation for the $k^\text{th}$ mutation at $v_{i;j}$ is
\eq{
A_{i;j} \mu^*_{i;j} \lrp{A_{i;j}'} = 
\begin{cases}
A_{i+1;j+1} A_{i-1;j-1} + A_{i;j+1} A_{i;j-1} & \mbox{if } i \neq 1 \mbox{ and } j \neq i+1\\
A_{i+1;j+1} A_{k;k} + A_{i;j+1} A_{i;j-1} & \mbox{if } i = 1 \mbox{ and } j \neq i+1\\
A_{i+1;j+1} A_{i-1;j-1} + A_{i;j+1} A_{k+i;k+i} & \mbox{if } i \neq 1 \mbox{ and } j = i+1\\
A_{i+1;j+1} A_{k;k} + A_{i;j+1} A_{k+i;k+i} & \mbox{if } i = 1 \mbox{ and } j = i+1.
\end{cases}
}{eq:mutij}
The variable $A_{i;j}$ is only replaced via a mutation at $v_{i;j}$,
so let's write $A_{{i;j}_k}$ to denote the unfrozen variable obtained after the $k^\text{th}$ mutation at $v_{i;j}$.
Then (\ref{eq:mutij}) becomes 
\eq{
A_{{i;j}_{k-1}} \mu^*_{i;j} \lrp{A_{{i;j}_k}} = 
\begin{cases}
A_{{i+1;j+1}_{k-1}} A_{{i-1;j-1}_k} + A_{{i;j+1}_{k-1}} A_{{i;j-1}_k} & \mbox{if } i \neq 1 \mbox{ and } j \neq i+1\\
A_{{i+1;j+1}_{k-1}} A_{k;k} + A_{{i;j+1}_{k-1}} A_{{i;j-1}_k} & \mbox{if } i = 1 \mbox{ and } j \neq i+1\\
A_{{i+1;j+1}_{k-1}} A_{{i-1;j-1}_k} + A_{{i;j+1}_{k-1}} A_{k+i;k+i} & \mbox{if } i \neq 1 \mbox{ and } j = i+1\\
A_{{i+1;j+1}_{k-1}} A_{k;k} + A_{{i;j+1}_{k-1}} A_{k+i;k+i} & \mbox{if } i = 1 \mbox{ and } j = i+1.
\end{cases}
}{eq:mutijk}
Recall that $A_{{i;j}_0} = \Delta^{1,\dots,i}_{j-i+1,\dots,j}$, and note that for $A_{{1;2}_1}$ this states
\eqn{
\Delta^1_2\, \mu^*_{1;2}\lrp{A_{{1;2}_1}} =  
 \Delta^1_1 \Delta^{1,2}_{2,3}+ \Delta^1_3 \Delta^{1,2}_{1,2}.
}
By Lemma~\ref{lem:minorIdent}, $\mu^*_{1;2}\lrp{A_{{1;2}_1}}= \Delta^{1,2}_{1,3}$.
Now use the mutation sequence to order the indices ${i;j}_k$.
Suppose the claim holds for all ${i';j'}_{k'} < {i;j}_k$.
As a notational convenience, let $\mu$ (without any subscripts) denote the birational map $T_{N;s_0} \dashrightarrow T_{N;s}$
from the torus of the initial seed to the torus of whichever seed $s$ is currently of interest, 
with $\mu$ obtained as the composition of a given sequence of cluster mutations. 
Then for $i\neq 1$ and $j\neq i+1$, $\Delta^{1,\dots,k+i-1}_{1,\dots,k-1,k+j-i,\dots,k+j-1}\, \mu^* \lrp{A_{{i;j}_k}} $ is given by
\eqn{
\Delta^{1,\dots,k+i-1}_{1,\dots,k,k+j-i+1,\dots,k+j-1} 
\Delta^{1,\dots,k+i}_{1,\dots,k-1,k+j-i,\dots,k+j} 
 +
\Delta^{1,\dots,k+i-1}_{1,\dots,k-1,k+j-i+1,\dots,k+j}
\Delta^{1,\dots,k+i}_{1,\dots,k,k+j-i,\dots,k+j-1}.
}
Using Lemma~\ref{lem:minorIdent} with $j_1 = k$, $j_2 = k+j-i$, and $j_3=k+j$, we find 
$$\mu^*\lrp{A_{{i;j}_k}}= \Delta^{1,\dots, k+i}_{1,\dots,k,k+j-i+1,\dots, k+j},$$
in agreement with the claim.
The remaining cases are essentially simplified versions of this case.
For $i = 1$ and $j\neq i+1$, $\Delta^{1,\dots,k}_{1,\dots,k-1,k+j-1}\, \mu^* \lrp{A_{{1;j}_k}} $ is given by
\eqn{
\Delta^{1,\dots,k}_{1,\dots,k} 
\Delta^{1,\dots,k+1}_{1,\dots,k-1,k+j-1, k+j} 
+
\Delta^{1,\dots,k}_{1,\dots,k-1,k+j}
\Delta^{1,\dots,k+1}_{1,\dots,k,k+j-1}. 
}
Taking $j_1 = k$, $j_2 = k+j-1$, and $j_3=k+j$,
Lemma~\ref{lem:minorIdent} yields
$$\mu^*\lrp{A_{{1;j}_k}}= \Delta^{1,\dots, k+1}_{1,\dots,k,k+j},$$
as desired.
For $i \neq 1$ and $j = i+1$, $\Delta^{1,\dots,k+i-1}_{1,\dots,k-1,k+1,\dots,k+i}\, \mu^* \lrp{A_{{i;i+1}_k}} $
is given by
\eqn{
\Delta^{1,\dots,k+i-1}_{1,\dots,k,k+2,\dots,k+i} 
\Delta^{1,\dots,k+i}_{1,\dots,k-1,k+1,\dots,k+i+1}
+
\Delta^{1,\dots,k+i-1}_{1,\dots,k-1,k+2,\dots,k+i+1}
\Delta^{1,\dots,k+i}_{1,\dots,k+i}.  
}
Now $j_1=k$, $j_2 = k+1$, and $j_3 = k+i+1$, and 
$$\mu^*\lrp{A_{{i;i+1}_k}}= \Delta^{1,\dots, k+i}_{1,\dots,k,k+2,\dots, k+i+1}.$$
Finally, for $i=1$ and $j=2$, $\Delta^{1,\dots,k}_{1,\dots,k-1,k+1}\, \mu^* \lrp{A_{{1;2}_k}} $ is given by
\eqn{
\Delta^{1,\dots,k}_{1,\dots,k} 
\Delta^{1,\dots,k+1}_{1,\dots,k-1,k+1,k+2} 
+ 
\Delta^{1,\dots,k}_{1,\dots,k-1,k+2}
\Delta^{1,\dots,k+1}_{1,\dots,k+1}. 
}
With $j_1 =k$, $j_2 = k+1$, and $j_3= k+2$, we find
$$\mu^*\lrp{A_{{1;2}_k}}= \Delta^{1,\dots, k+1}_{1,\dots,k,k+2}.$$
By induction, Proposition~\ref{prop:Gew0minors} holds.
\end{proof}

{\cor{
In $\mathring{U}$,
$A_{{i;j}_k}$ pulls back to $T_{N;s_0}$ as 
$\Delta^{k+1,\dots, k+i}_{k+j-i+1,\dots, k+ j}$.
In particular, for the final seed $s'$ reached by the mutation sequence,
$\mu^*\lrp{A_{{i;j}_{s'}}} = \alpha^*\lrp{A_{{i;j}_{s_0}}}$.
}}

\begin{proof}
The variables on $\mathring{U}$ are restrictions of the variables on $G^{e,w_0}$.
On $U$, $\Delta^{1,\dots,k,k+1, \dots, k+i}_{1,\dots,k,k+j-i+1,\dots,k+j} = \Delta^{k+1, \dots, k+i}_{k+j-i+1,\dots,k+j}$,
so the first part follows immediately from Proposition~\ref{prop:Gew0minors}. 
Then every mutation at $v_{i;j}$ in the sequence adds 1 to every row and column index of the minor $\mu^*\lrp{A_{i;j}}$.
There are $n-j$ mutations at $v_{i;j}$, the row indexing set $I=\lrc{1,\dots, i}$ for $A_{{i;j}_{s_0}}$ becomes 
\eqn{
\lrc{n-j+1, \dots, n-j+i} = 
\lrc{n+1-j, \dots, n+1-(j-i+1)} =
\lrc{w_0(j), \dots, w_0(j-i+1)}
}
Similarly, the column indexing set $J = \lrc{j-i+1, \dots, j}$ becomes
\eqn{
\lrc{n-j+j-i+1, \dots, n-j+j} = 
\lrc{n+1-i, \dots, n+1-1} =
\lrc{w_0(i), \dots, w_0(1)}.
}
So $I$ is sent $w_0(J)$ with the opposite ordering of indices and $J$ is sent to $w_0(I)$ with the opposite ordering of indices.
Since both orders have been reversed, no sign is introduced.
Then 
\eqn{
\mu^*\lrp{A_{{i;j}_{s'}}} = \Delta^{w_0(J)}_{w_0(I)} = \alpha^*\lrp{A_{{i;j}_{s_0}}}.
} 
\end{proof}

This completes the proof of Proposition~\ref{prop:2torimutequiv}.

\subsection{Full Fock-Goncharov conjecture for $U$}
Many sufficient conditions for the full Fock-Goncharov conjecture to hold are provided in \cite{GHKK}
For $U$, they spelled things out very explicitly in Conjecture 11.11 and Proposition 11.12.

{\bf{\cite{GHKK} Conjecture 11.11.}}  {\it{Each of the (non-constant) matrix entries is a cluster variable for
some cluster of U. There is a choice of seed for which the $\gv$-vectors for the matrix entries are linearly independent.
}}

{\bf{\cite{GHKK} Proposition 11.12.}} {\it{\cite{GHKK} Conjecture 11.11 implies the full Fock-Goncharov conjecture for $U$.}}

I will establish \cite{GHKK} Conjecture 11.11 here.

{\prop{\label{prop:vars}
Every nonconstant matrix entry is a cluster variable in some seed.}}

\begin{proof}
In $U$, $\Delta^{1,\dots,i}_{1,\dots,i-1,j}$ represents the determinant of an upper triangular matrix with the property that only the last diagonal entry may differ from $1$. 
Then the following simplification holds in $U$:
\eqn{
\Delta^{1,\dots,i}_{1,\dots,i-1,j} = \Delta^{i}_{j}. 
}
With this in mind, it is sufficient to show that for $i<j$, 
$\Delta^{1,\dots,i}_{1,\dots,i-1,j}$ 
is a cluster variable of $G^{e,w_0}$.
$\SL_n$ is of Cartan type $A_{n-1}$, so it is simply laced. 
Remark 2.14 of \cite{FZ_clustersIII} indicates that every pair of seeds in $G^{e,w_0}$ 
coming from reduced words for $w_0$ are mutation equivalent.
It suffices to produces a word that yeilds 
$\Delta^{1,\dots,i}_{1,\dots,i-1,j}$ 
as a chamber minor.
Any word resulting in some region where the lowest $i$ blue lines are 
\textcolor{blue}{$1_B$}, \dots, 
\textcolor{blue}{$(i-1)_B$},
and
\textcolor{blue}{$j_B$}
would produce the desired chamber.
I proceed by induction.
There are no unfrozen variables for $n<3$, so take $n\geq 3$.
For $n=3$, only $\Delta^{2}_{3}$ is not in the initial seed.
However, this is a chamber minor for the reduced word $\lrp{2,1,2}$.
Now suppose the result holds for all $m<n$.

If $j_B<n$, start the reduced word with the sequence $1,2,\dots, n-1$.
This moves $n_B$ to its final position at the top of the double pseudoline arrangement.
The problem is reduced to the case of $n-1$, where the induction hypothesis applies.
Finally, if $j_B = n$, start the word with $2,3,\dots,n-1,2,3,\dots,n-2,\dots,n-i+1, n-i+2, \dots, i$.
This moves $(n-1)_B$ to the top, then $(n-2)_B$ below it, and so on until $i_B$.
The lines are in the desired position.
To finish the word, bring $n_B$ to the top via $1,2, \dots, n-1$.
Now $i_B$ through $n_B$ are in their final positions, and no pair of lines has crossed more than once.
Appending any reduced word for the longest element in the $n=i-1$ case completes the job. 
\end{proof}

Alternatively, recall that in the sequence of mutations described in Proposition~\ref{prop:2torimutequiv},
each mutation at $v_{i;j}$ increased every row index and column index for the minor $\mu^*\lrp{A_{i;j}}$ by 1.
Then for $1\leq k < l \leq n$, 
$\Delta^k_l$ is obtained as $\mu^*\lrp{A_{{1;k-l+1}_{k-1}}}$.
This explicit construction will be useful for computing the $\gv$-vectors of the non-constant matrix entries.
Before doing so, I provide a description of $\gv$-vectors for the unfamiliar reader.
The description provided here is intended only as an overview, and is taken from 
\cite{GHK_birational} and \cite{GHKK}.
$\gv$-vectors were introduced in 
\cite{FZ_clustersIV}
with a flavor quite different from Gross, Hacking, Keel, and Kontsevich's description,
but {\it{almost}} equivalent content.
They use a different convention in their mutations of ``cluster variables with principal coefficients'' 
and as a result, explicit calculations of $\gv$-vectors within the two frameworks will differ.\footnote{This difference should be noted before wagering a pitcher of beer on the correct $\gv$-vectors of cluster monomials.}

Let $\lrp{\cA,\cX}$ be the pair of cluster varieties in a cluster ensemble, and for simplicity take all multipliers $d_i$ of the inital data to be 1.
There is a larger space, denoted $\cAp$, that (roughly speaking) is built out of this pair.
If $N$ is the cocharacter lattice for (a torus in the atlas of) $\cA$, 
its dual lattice $M$ plays the corresponding role for $\cX$.
The space $\cAp$ is the $\cA$-variety with cocharacter lattice $\widetilde{N}:= N\oplus M$ and
skew form 
\eqn{
\lrc{ \lrp{n_1,m_1}, \lrp{n_2,m_2} }_{\widetilde{N}} := \lrc{n_1,n_2}_N + \lra{n_1,m_2} - \lra{n_2,m_1}.
}
If $s_0 = \lrc{e_1, \dots, e_n}$ is the initial seed of $\cA$, 
$\widetilde{s_0}= \lrc{ \lrp{e_1, 0},\dots, \lrp{e_n, 0}, \lrp{0, e_1^*},\dots, \lrp{0, e_n^*}}  $
is the initial seed for $\cAp$.
The unfrozen sublattice $\widetilde{N}_{\mathrm{uf}}$ is identified with the unfrozen sublattice $N_{\mathrm{uf}}$ of $N$.
A choice of seed $s$ determines a canonical partial compactification $\cAps{s}$ of $\cAp$ 
by allowing the frozen 
$X_{i_s} := z^{\lrp{0,e_i^*}_s}$
variables to vanish.
There is a canonical surjection $\pi: \cAps{s} \to \C^n_{X_{1_s}, \dots, X_{n_s}}$,
defined on cocharacter lattices by $\lrp{n,m} \mapsto m$.
$\cA$ is identified with $\pi^{-1} \lrp{1,1, \dots,1} \subset \cAps{s}$.
The choice of seed $s$ also yields a canonical extension of cluster monomials\footnote{That is, regular functions on $\cA$ that are monomials on some torus in the atlas for $\cA$.}
from $\cA$ to $\cAp$.
$z^m$ extends to $z^{\lrp{m,0}}$, and if the function is a monomial on $T_{N;s'}$, 
mutation from $s'$ to $s$ can be performed equally well in $\cA$ and $\cAp$, yielding the extension of the cluster monomial.
The Laurent phenomenon indicates that the cluster monomial further extends to $\cAps{s}$.
The map 
\eqn{
N &\to N \oplus M\\
n &\mapsto \lrp{n,p^*(n)}.
} 
on the level of cocharacter lattices induces an action of $T_N$ on $\cAps{s}$.
The extension to $\cAps{s}$ of a cluster monomial is an eigenfunction of this action.
The {\it{$\gv$-vector at $s$}} of the cluster monomial is its weight under this action, identified with some point in $M$.
Alternatively (and more simply in practice), the extension of the cluster monomial restricts to a regular non-vanishing eigenfunction along the central fiber $\pi^{-1}(0)$.
The $\gv$-vector is the $T_N$-weight of this restriction to the central fiber.
Essentially, this alternative description allows us to drop all terms in the 
cluster monomial having some $X_{i_s}$ factor and compute the $\gv$-vector of 
the more manageable resulting expression.

{\prop \label{prop:g-vects}
{Let $s$ be the seed 
corresponding to the image of $\varphi_{{\bf{i}}_{w_0}}$.
The ${\gv}$-vector at $s$ for the minor $\Delta^{1,\dots,k+i}_{1,\dots,k,k+j-i+1,\dots,k+j}= \mu^*\lrp{A_{{i;j}_k}}$ in $\ssO\lrp{G^{e,w_0}}$ is 
$e^*_{n+i-j-k;n-k} - e^*_{n-j-k;n-k-i} + e^*_{i+k;i+k}$, 
where $e^*_{r;s}$ is interpreted as $0$ if $v_{r;s}$ is not a vertex of the quiver for $G^{e,w_0}$.
Furthermore, setting $e^*_{r;r} = 0$ yields the $\gv$-vector for the minor in $\ssO(U)$.
In particular, the $\gv$-vector at $s$ for the
non-constant matrix entry $\alpha^*\lrp{\Delta^{i}_{j}} = \Delta^{n+1-j}_{n+1-i}$ is $e^*_{i;j}$ if $i=1$ and $e^*_{i;j} - e^*_{i-1;j-1}$ otherwise.
}}

\begin{proof}
This minor is $\mu^*\lrp{A_{{i;j}_k}}$.
The quiver for $s$ is obtained from the quiver in Figure~\ref{fig:Gew0_quiver_final} (upon relabeling the vertices)
by introducing a frozen vertex $w_{i;j}$ for each unfrozen $v_{i;j}$, and an arrow $v_{i;j} \to w_{i;j}$.
Denote by $W$ the set of these new frozen vertices $w_{i;j}$.
I show the resulting quiver for $n=4$ below.

\input{Gew0FinalDouble.tex}

The quiver at every other seed in the mutation pattern can be obtained by mutating this quiver through the inverse mutation pattern.

{\lemma{\label{lem:vWarrows}
Let $Q_{\text{prin},s}$ be the quiver for the seed $s$ of $G^{e,w_0}$ with principal coefficients at $s$,
and let $\nu$ be the inverse mutation pattern to that described in Proposition~\ref{prop:2torimutequiv}.
Then the quiver after the $k^{\text{th}}$ mutation at $v_{i;j}$ in $\nu$ 
has exactly the following arrows between $v_{i;j}$ and $W$:
$\lrc{w_{i+k-r;j+k-r} \to v_{i;j} : 1\leq r \leq i}$.
}}

\begin{proof}
Add an index to count mutations at each vertex, so ${v_{i;j}}_k$ corresponds to the $k^{\text{th}}$ mutation at $v_{i;j}$.
Then use $\nu$ to order the set $\lrc{{v_{i;j}}_k}$: ${v_{1;2}}_1 < {v_{2;3}}_1 < {v_{1;3}}_1 \cdots$.
The claim clear for ${v_{1;2}}_1$.
Assume it holds for all ${v_{i';j'}}_{k'} < {v_{i;j}}_k$.
Only mutations at $v_{i;j}$ or an adjacent vertex may introduce changes to the set of arrows for which $v_{i;j}$ is a vertex (and mutations at $v_{i;j}$ just change the direction of those arrows).
Recall the mutations of Lemma~\ref{lem:delbot} and Corollary~\ref{cor:revarrows}.
The mutation ${v_{i';j'}}_{k'}$ is at a vertex adjacent to $v_{i;j}$ if and only if 
$v_{i';j'} \in \lrc{ v_{i+1;j+1},v_{i-1;j-1}, v_{i;j+1}, v_{i;j-1}}$.
Moreover, the arrow connecting $v_{i';j'}$ and $v_{i;j}$ is directed away from $v_{i;j}$ 
in the first two case and toward $v_{i;j}$ in the last two. 
I'll treat the $k>1$ case first.
By the induction hypothesis, after ${v_{i;j}}_{k-1}$, the arrows between $v_{i;j}$ and $W$ are
$\lrc{w_{i+k-1-r;j+k-1-r} \to v_{i;j} : 1\leq r \leq i}$.
Prior to ${v_{i+1;j+1}}_{k-1}$, the arrows between $v_{i+1;j+1}$ and $W$ are
$\lrc{v_{i+1;j+1} \to w_{i+k-r;j+k-r}: 1\leq r \leq i+1}$.
Before cancelling 2-cycles, composition yields arrows 
$\lrc{v_{i;j} \to w_{i+k-r;j+k-r}: 1\leq r \leq i+1}$.
All that remains upon cancelling 2-cycles is
$v_{i;j} \to w_{i+k-1;j+k-1}$.
No compositions yielding arrows between $v_{i;j}$ and $W$ can occur from mutation at $v_{i;j+1}$ or $v_{i;j-1}$
as the arrows are directed from these vertices and directed toward $v_{i;j}$ and $W$.
More compositions do arise from $v_{i-1;j-1}$, if $i>1$.
Prior to ${v_{i-1;j-1}}_k$, the arrows between $v_{i-1;j-1}$ and $W$ are
$\lrc{v_{i-1;j-1} \to w_{i-1+k-r;j-1+k-r}: 1\leq r \leq i-1}$.
After performing this mutation, 
the set of arrows between $v_{i;j}$ and $W$ becomes
$\lrc{v_{i;j} \to w_{i-1+k-r;j-1+k-r} : 0 \leq r \leq i-1} 
=
\lrc{v_{i;j} \to w_{i+k-r;j+k-r} : 1 \leq r \leq i}
$.
Finally, the mutation ${v_{i;j}}_k$ reverses each of these arrows, yielding the claim for $k>1$, assuming the $k=1$ result.
A simplified version of the $k>1$ argument applies for $k=1$.
Now no mutation has occurred at $v_{i+1;j+1}$, so we skip to mutation at $v_{i-1;j-1}$, if $i>1$.
Then prior to mutation at $v_{i;j}$,
the set of arrows between $v_{i;j}$ and $W$ is
$\lrc{v_{i;j} \to w_{i+1-r;j+1-r} : 1 \leq r \leq i}$.
Mutation at $v_{i;j}$ reverses these.
\end{proof}

The variables at $s$ are monomials on $T_{N;s}$: $A_{{i;j}_s}$ is just $z^{e^*_{i;j}}$ on this torus.
So $\gv_s\lrp{A_{{i;j}_s}} = e^*_{i;j}$.
The cluster variables from other seeds need to be pulled back to $s$ via $\mu^{-1}$.
The terms on the right of (\ref{eq:mutij}) are unaltered by $\mu_{i;j}$,
so ${\mu_{i;j}^{-1}}^*\lrp{A_{i;j}} $ can be read off of (\ref{eq:mutij}).
Upon including the vertices for $X$ variables, the second term in each binomial on the right of (\ref{eq:mutij}) gains factors of $X$ variables, while the first term in the binomial does not.
So over the central fiber, 
\eqn{
{\mu_{i;j}^{-1}}^*\lrp{A_{i;j}} A_{i;j}'=  \begin{cases}
A_{i+1;j+1}' A_{i-1;j-1}' & \mbox{ if } i\neq 1\\
A_{i+1;j+1}' A_{k;k}' & \mbox{ if } i = 1. 
\end{cases}
}
Then if $\mu_{{i;j}_k}: s' \to s''$ is the $k^\text{th}$ mutation at $v_{i;j}$,
\eqn{
\gv_s\lrp{A_{{i;j}_{s'}}} = \gv_s\lrp{A_{{i;j}_{k-1}}} &= \begin{cases}
\gv_s \lrp{A_{{i+1;j+1}_{s''}}}
+\gv_s \lrp{A_{{i-1;j-1}_{s''}}}
-\gv_s \lrp{A_{{i;j}_{s''}}}
& \mbox{ if } i\neq 1\\
\gv_s \lrp{A_{{i+1;j+1}_{s''}}}
+\gv_s \lrp{A_{k;k}}
-\gv_s \lrp{A_{{i;j}_{s''}}}
& \mbox{ if } i= 1
\end{cases}\\
& = \begin{cases}
\gv_s \lrp{A_{{i+1;j+1}_{k-1}}}
+\gv_s \lrp{A_{{i-1;j-1}_{k}}}
-\gv_s \lrp{A_{{i;j}_{k}}}
& \mbox{ if } i\neq 1\\
\gv_s \lrp{A_{{i+1;j+1}_{k-1}}}
+\gv_s \lrp{A_{k;k}}
-\gv_s \lrp{A_{{i;j}_{k}}}
& \mbox{ if } i= 1.
\end{cases}
}

The variables of the seed $s$ are
$A_{{i;j}_{w_0(j)-1}}$.
Then
\eqn{
\gv_s\lrp{A_{{i;j}_{w_0(j)-1}}} = e^*_{i;j}.
}
This matches the claim:
\eqn{
e^*_{n+i-j-\lrp{w_0(j)-1};n-\lrp{w_0(j)-1}} 
- e^*_{n-j-\lrp{w_0(j)-1};n-\lrp{w_0(j)-1}-i}
&= 
e^*_{n+i-j-\lrp{n-j};n-\lrp{n-j}} 
- e^*_{n-j-\lrp{n-j};n-\lrp{n-j}-i}\\
&= 
e^*_{i;j} 
- e^*_{0;-i}\\
&=
e^*_{i;j} 
}
since $v_{0;-i}$ is not a vertex in the quiver for $G^{e,w_0}$.
Now use the sequence of mutations to order the indices ${i;j}_k$ so that 
${i;j}_k < {i';j'}_{k'}$ if the 
$k^{\text{th}}$ mutation at $v_{i;j}$ occurs 
before the 
$k'^{\text{ th}}$ mutation at $v_{i';j'}$, and set the indices for all frozen variables to be greater than the indices of unfrozen variables.
Suppose the claim holds for all ${i';j'}_{k'} > {i;j}_{k-1}$.
Then $\gv_s\lrp{A_{{i;j}_{k-1}}}$ is given by  
\eqn{
& e^*_{n+(i+1)-(j+1)-(k-1); n-(k-1)} -e^*_{n-(j+1)-(k-1);n-(k-1)-(i+1)}  + e^*_{(i+1) + (k-1);(i+1)+(k-1)}\\
&+e^*_{n+(i-1)-(j-1)-k; n-k} -e^*_{n-(j-1)-k;n-k-(i-1)}  +e^*_{(i-1)+k;(i-1) +k}\\
&-e^*_{n+i-j-k; n-k} +e^*_{n-j-k;n-k-i}  -e^*_{i+k;i+k}\\
= 
& e^*_{n+i-j-(k-1); n-(k-1)} -e^*_{n-j-k;n-k-i} +e^*_{i+k;i+k}  \\
&+e^*_{n+i-j-k; n-k} -e^*_{n-j-(k-1);n-(k-1)-i}+e^*_{i+(k-1);i +(k-1)}\\
&-e^*_{n+i-j-k; n-k} +e^*_{n-j-k;n-k-i}-e^*_{i+k;i+k}\\
=
& e^*_{n+i-j-(k-1); n-(k-1)} 
-e^*_{n-j-(k-1);n-(k-1)-i}+e^*_{i+(k-1);i +(k-1)}
}
if $i\neq 1$, 
and 
\eqn{
& e^*_{n+2-(j+1)-(k-1); n-(k-1)} -e^*_{n-(j+1)-(k-1);n-(k-1)-2}+ e^*_{2 + (k-1);2+(k-1)}\\
&+e^*_{k;k}\\
&-e^*_{n+1-j-k; n-k} +e^*_{n-j-k;n-k-1}-e^*_{1+k;1+k}\\
= 
& e^*_{n+1-j-(k-1); n-(k-1)} -e^*_{n-j-k;n-k-1}+ e^*_{1 + k ;1 + k}\\
&+e^*_{k;k}\\
&-e^*_{n-j-(k-1); n-k} +e^*_{n-j-k;n-k-1}-e^*_{1+k;1+k}\\
=
& e^*_{n+1-j-(k-1); n-(k-1)} 
-e^*_{n-j-(k-1); n-(k-1)+1} +e^*_{k;k}
}
if $i=1$.
This proves the first claim.  The remaining claims follow straightfowardly from the first. 
\end{proof}

{\cor \label{cor:FG_U}
{The full Fock-Goncharov conjecture holds for $U$.}}

\begin{proof}
For the non-constant matrix entries, 
the $\gv$-vectors at $s$ described in Proposition~\ref{prop:g-vects} are linearly independent.
Together with Proposition~\ref{prop:vars}, 
this establishes \cite{GHKK} Conjecture 11.11.
Then \cite{GHKK} Proposition 11.12 yields the claim.
\end{proof}

\subsection{Potential and canonical basis for $U$}
$U$ is obtained from $\mathring{U}$ by allowing the frozen variables $A_{i;n}=\Delta^{1,\dots, i}_{n-i+1,\dots, n}$ to vanish.
Corresponding to each divisor $A_{i;n} = 0$ is a cluster monomial on the dual $\cX$ variety.
If $v_{i;n}$ is a sink in the quiver $Q_{s'}$ of some seed $s'$,
then $z^{-e_{{i;n}_{s'}}}$ pulls back to a regular function on every torus adjacent to $T_{M;s'}$.
Up to a codimension 2 subset, this is all of $\cX$, so 
$z^{-e_{{i;n}_{s'}}}$
would extend to a global regular function on $\cX$.\footnote{%
See \cite{GHKK} Section 10 for a more complete discussion.
}
This function is denoted $\vartheta_{i;n}$, and a seed where $v_{i;n}$ is a sink is said to optimized for ${i;n}$.

The $\vartheta$-functions for frozen indices are key to describing the canonical basis for $U$.
A Landau-Ginzburg potential $W :\cX \to \C$ is constructed as the sum of these $\vartheta$-functions, 
and the canonical basis of $H^0\lrp{U,\ssO_U}$ is parametrized by $W^T\geq 0 \subset \cX\lrp{\Z^T}$.\cite[Corollary~11.9,Corollary~11.10]{GHKK}
The non-negative locus of $\vartheta_{i;n}^T$ 
yields the subset of $\cX\lrp{\Z^T}$ that is regular along the divisor $A_{i;n}=0$.
Background on tropicalizations is described cleanly in \cite[Section~2]{GHKK}. 

{\prop{\label{prop:opt_seeds}
Every frozen index for $G^{e,w_0}$ has an optimized seed.
As a result, each frozen index for $U$ has an optimized seed.}}

\begin{proof}
Consider a quiver $Q_L$ of the form
\input{QLine.tex}
The sequence of mutations $v_1,v_2, \dots, v_n$ yields the quiver
\input{QLineMut.tex}
making $v_0$ a sink.
The quiver at $s$ for $G^{e,w_0}$ is shown in Figure~\ref{fig:Gew0_quiver_final}, up to the previously described relabeling.
Call it $Q_{G^{e,w_0}}$.
This quiver is optimized for $v_{n;1}$ and $v_{n-1;n}$.
For the remaining frozen vertices $v_f$, there is a subquiver of $Q_{G^{e,w_0}}$ isomorphic to $Q_L$ in which $v_f$ plays the role of $v_0$.
Performing these mutations on $Q_{G^{e,w_0}}$ only affects the subquiver whose vertices are either in $Q_L$ or connected to $Q_L$ by an arrow.
As arrows between frozen vertices are deleted,
any frozen vertices besides $v_f$ can be ignored when determining if $v_f$ becomes a sink.
Then the relevant subquiver of $Q_{G^{e,w_0}}$ has the form 
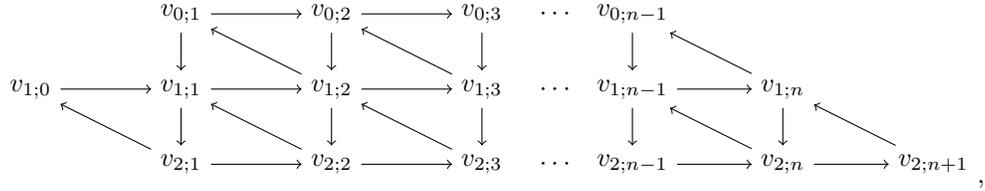
\input{QSubFull.tex}
possibly with the top or bottom row deleted and with the middle row being the subquiver $Q_L$.
The cycles prevent any new arrows involving $v_{1;0}$ from developing via some composition with an arrow not in $Q_L$.
The explicit mutations are shown below.
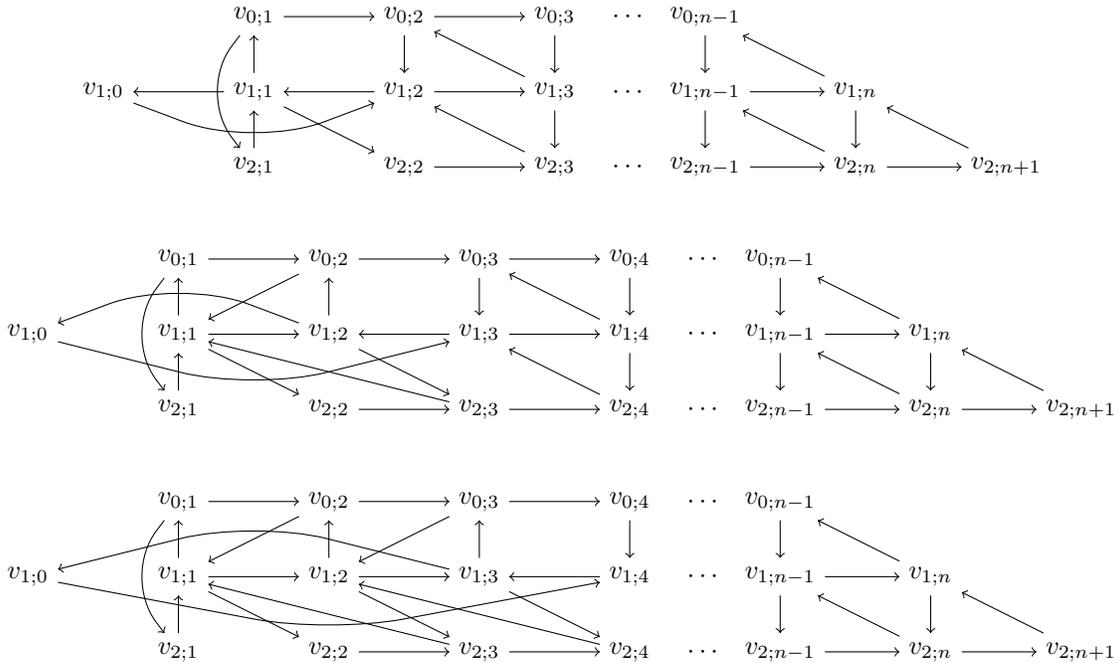
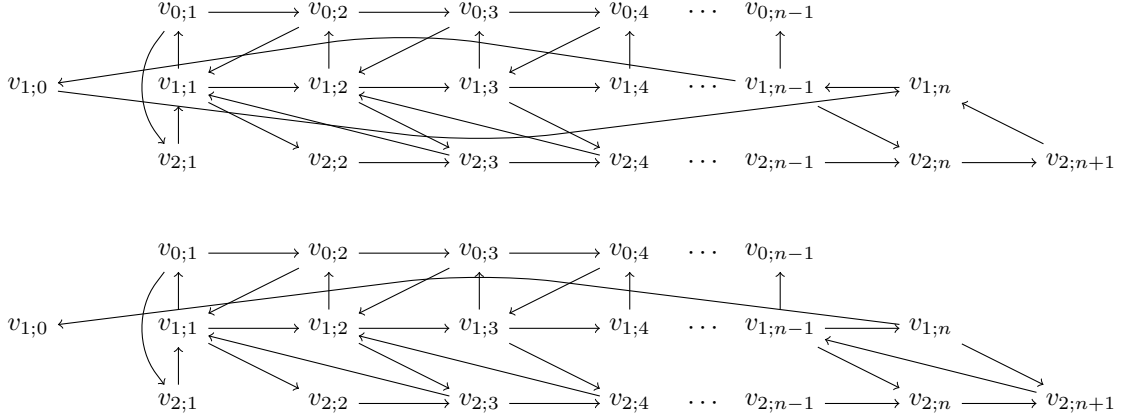
\input{QSubFullmut1.tex}
\input{QSubFullmut2.tex}
\input{QSubFullmut3.tex}
\eqn{\vdots}
\input{QSubFullmutn-1.tex}
\input{QSubFullmutn.tex}
\end{proof} 

If $s_{\text{op}}$ is the seed optimized for $v_f$ described above, 
then $\vartheta_f$ restricts to $T_{M;s_{\text{op}}}$ as $z^{-e_f}$.
Pulling back via the sequence of birational maps associated to the quiver mutations yields
the restriction of $\vartheta_f$ to $T_{M;s}$.
The formula for mutation at $v_k$ is 
\eqn{
\mu_k^*\lrp{ z^n } = z^n \lrp{ 1 + z^{e_k}}^{-\lrc{n,e_k}},  
}
and if $s_1= \lrc{e_1, \dots, e_n}$, then $\mu_k \lrp{s_1} = \lrc{e_1',\dots, e_n'}$ where
\eqn{
e_i' = \begin{cases}
e_i + \lrb{\epsilon_{ik}}_+ e_k & {\mbox{if }} i\neq k\\
-e_k & {\mbox{if }} i=k
\end{cases}.
}

Note that the pullback of 
$z^{-e_{f_{s_{\text{op}}}}}$
to $T_{M;s}$ 
will involve only the subset of $s$ corresponding to vertices at which a mutation has occurred and the vertex $v_f$.
So only the subquiver $Q_L$ is needed to determine the pullback.

{\prop{\label{prop:theta}
Using the indexing for $Q_L$ at the beginning of Proposition~\ref{prop:opt_seeds},
the pullback of $z^{-e_{0_{s_{\text{op}}}}}$ to $T_{M;s}$ is 
\eqn{
z^{-e_0} + 
z^{-e_0-e_1} + 
z^{-e_0-e_1-e_2} + 
\cdots +
z^{-e_0-e_1-e_2 -\cdots - e_n}. 
}
}}

\begin{proof}
The quiver for the first mutation is 
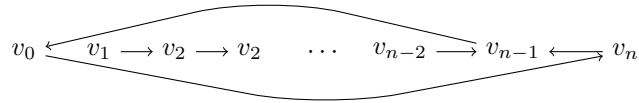
\input{QLineMut1.tex}
So 
\eqn{
\mu_{n}^* \lrp{z^{-e_0'}} &= z^{-e_0'}\lrp{1+z^{e_n}}^{-\lrc{-e_0',e_n}}\\
&= z^{-e_0 - \lrb{ \epsilon_{0,n}}_+ e_n} \lrp{1 + z^{e_n}}^{\lrc{e_0 + \lrb{\epsilon_{0,n}}_+e_n,e_n}}\\
&= z^{-e_0 - e_n} \lrp{1 + z^{e_n}}.
}
The next quiver is 
\input{QLineMut2.tex}
\eqn{
\mu_{n-1}^* \lrp{z^{-e_0' - e_n'} \lrp{1 + z^{e_n'}}} &= 
z^{-e_0' - e_n'}\lrp{1+z^{e_{n-1}}}^{-\lrc{-e_0'-e_n',e_{n-1}}} \lrp{1 + z^{e_n'} \lrp{1+z^{e_{n-1}}}^{-\lrc{e_n',e_{n-1}}}}\\
&=
z^{-e_0 - e_{n-1} - e_n }\lrp{1+z^{e_{n-1}}}^0 \lrp{1 + z^{e_n} \lrp{1+z^{e_{n-1}}}^{1}}\\
&=
z^{-e_0 - e_{n-1} - e_n }\lrp{1 + z^{e_n} \lrp{1+z^{e_{n-1}}}}.
}
This pattern continues with the $i^{\text{th}}$ mutation yielding
\eqn{
z^{-e_0-e_{n-i+1} - e_{n-i+2} - \cdots - e_{n}} \lrp{1 + z^{e_n} \lrp{ 1 + z^{e_{n-1}} \cdots \lrp{1+z^{e_{n-i+1}}} \cdots }}. 
}
The result after all $n$ mutations is
\eqn{
& z^{-e_0-e_1 - \cdots - e_n} \lrp{1 + z^{e_n} \lrp{ 1 + z^{e_{n-1}} \cdots \lrp{1+z^{e_{1}}} \cdots }}\\
=& z^{-e_0-e_1 - \cdots - e_n}
+  z^{-e_0-e_1 - \cdots - e_{n-1}}
+ \cdots + z^{-e_0}, 
}
as claimed.
\end{proof}

Using Proposition~\ref{prop:theta}, the restriction of the potential $W$ to $T_{M;s}$ can immediately be written down for $U$.

{\cor{\label{cor:W}
The restriction of $W$ to $T_{M;s}$ is
\eqn{
\sum_{i=1}^{n-1} \vartheta_{i;n},
}
where 
\eqn{
\vartheta_{i;n} = \sum_{j=0}^{i-1} z^{-\sum_{k=0}^j e_{i-k;n-k}}.
}
}
} 

The canonical basis for $U$ is obtained by tropicalizing this potential and taking its non-negative locus.
In the seed $s$, the basis has a very satisfying description.
{\prop{\label{prop:CanAndCone}
Let $\cX$ be the Fock-Goncharov dual of the $\cA$-variety $\mathring{U}$.
Then at $s$,
$$\Xi := \lrc{x \in \cX\lrp{\R^T} : W^T (x) \geq 0 } $$ 
is a full dimensional simplicial cone.
Furthermore, $\Xi$ is the $\R_{\geq 0}$ span of the $\gv$-vectors at $s$ for the non-constant matrix entries.
These $\gv$-vectors form a basis for $M$ (identified with $ \cX\lrp{\Z^T}$ by the choice of seed),
so they are exactly the non-identity generators of the monoid $\Xi \cap M$.
The canonical basis for $H^0\lrp{U,\ssO_U }$ is simply
\eqn{
B_U = \lrc{\prod_{{i;j} \in I} \lrp{\Delta^i_j}^{a_{i;j}} : a_{i;j} \in \Z_{\geq 0}}.
}  
}}

\begin{proof}
Recall the $\gv$-vectors at $s$ for the non-constant matrix entries
described in Proposition~\ref{prop:g-vects}.
The $\R_{\geq 0}$ span of this collection of vectors is clearly a full dimensional simplicial cone in $M_\R$.
So the second claim implies the first.
$W^T \geq 0$ defines a system of linear inequalities-- if $W = \sum_i c_i z^{n_i}$, $c_i \in \Z_{>0}$, then
\eqn{
W^T (x) \geq 0  \iff \min\lra{n_i , - x } \geq 0.  
}
I claim that $\lrc{-n_i}_i$ is the dual basis to the set of $\gv$-vectors in question.
First note that all $n_i$ with a non-zero $e_{j;n}$ component come from the same $\vartheta$-function,
and any $n_k$ with a non-zero $e_{j-l;n-l}$ component comes from this $\vartheta$-function as well.
The $\gv$-vectors in question are
\eqn{
\gv_s \lrp{\Delta^{w_0(j)}_{w_0(i)}} = \begin{cases}
e_{i;j}^* - e_{i-1;j-1}^* & \mbox{if } i\neq 1\\
e_{1;j}^* & \mbox{if } i=1
\end{cases}.
}
So the $\gv$-vectors are orthogonal to the subspaces spanned by the $n_i$ from all but one $\vartheta$-function.
It is sufficient to show that it pairs to $1$ with exactly one of the $n_i$ from this $\vartheta$-function, 
and it pairs to $0$ with the rest. 
For $e_{i;j}^*$, the relevant $\vartheta$-function is $\vartheta_{n+i-j;n}$.
\eqn{
\lra{
\sum_{k=0}^l e_{(n+i-j)-k;n-k} , 
e_{i;j}^* - e_{i-1;j-1}^*}
= \begin{cases}
1 & \mbox{if } l = n-j \\
0 & \mbox{otherwise}
\end{cases}
}
\eqn{
\lra{
\sum_{k=0}^l e_{(n+i-j)-k;n-k} , 
e_{1;j}^* }
= \begin{cases}
1 & \mbox{if } l =n-j \\
0 & \mbox{otherwise}
\end{cases}
}
The claim holds.
The rest of the proposition follows immediately.
\end{proof}
%
%
%

\section{Base affine space $G/U$}

In the previous section, the cluster variety structure of $U$ came through consideration of the double Bruhat cell $G^{e,w_0}$.
In that case, specializing a subset of the frozen chamber minors to 1 produced an open subset of $U$.
$G/U$ also inherits a cluster variety structure from $G^{e,w_0}$.
There is an open embedding of $G^{e,w_0}$ into $G/U$ given by
\eqn{
G^{e,w_0} &\hookrightarrow G/U\\
g &\mapsto g^T U.
}
This time no specialization is needed.

{\prop{\label{prop:Gew0InGmodU}
$g,g' \in G$ are in the same equivalence class of $G/U$ if and only if $g^T$ and $g'^T$ agree on all of the cluster variables of $G^{e,w_0}$.\footnote{%
To extend the variables to $G$ from the subset $G^{e,w_0}$, 
note that the variables can be written in terms of minors, which are defined on all of $G$.}
}}

\begin{proof}
Suppose $g =  g' u $ for some $u \in U$.
The minors parametrizing the initial seed are all determinants of top-aligned submatrices, that is, submatrices whose row indexing set has the form $1, 2, \dots, i$.
Say the submatrix in question is $r\times r$.
Since $u$ is upper triangular, only the top left $r \times r$ block of $u$ affects the relevant submatrix of $g^T=u^T g'^T$.
This block is again upper triangular with 1's on the diagonal.
So the submatrix of $g'^T$ has been multiplied by a determinant 1 submatrix of $u$.
$g^T$ and $g'^T$ agree on this minor, and on all minors in the initial seed.
Furthermore, all cluster variables are rational functions in the initial seed minors, 
so $g^T$ and $g'^T$ agree on all cluster variables.

Now suppose $g^T$ and $g'^T$ agree on all cluster variables.
$U$ is the stabilizer of the standard flag decorated with volume forms on each subspace in the flag.
The action of $g$ on the 1 dimensional decorated subspace of the decorated standard flag 
is determined by applying $g$ to each of the standard basis vectors. 
This corresponds to the determining the collection of $1\times 1$ left-aligned minors of $g$, or equivalently, top-aligned minors of $g^T$.
For the 2 dimensional decorated subspace, we need to see what happens to each parallelogram defined by a pair of basis vectors.
This corresponds to determining the collection of $2 \times 2$ top-aligned minors of $g^T$, and so forth.
For the $k$-th decorated subspace, we need the collection of top-aligned $k\times k$ minors.
If $g^T$ and $g'^T$ agree on each of these minors, then the image of the decorated standard flag is the same for the linear maps $g$ and $g'$,
so they are in the same equivalence class in $G/U$.
So it suffices to show that $\Delta^{1,\dots,k}_{j_1,\dots,j_k}$ is a cluster variable
for each $\lrc{j_1,\dots, j_k} \subset \lrc{1, \dots, n}$.
In fact, each such minor is a chamber minor.
We simply need to choose a reduced expression for $w_0$ where 
$\lrc{{j_1}_B,\dots, {j_k}_B} $ are the bottom blue lines for some portion of the double pseudoline arrangement.
The argument used to show that the minor in Proposition~\ref{prop:vars} is a chamber minor applies here as well.
The claim is proved.
\end{proof}

Now, $G^{e,w_0}$ is precisely the non-vanishing locus of certain minors in $G$, the frozen minors for the cluster variety.
So if we allow these minors to vanish, we get by the previous proposition a space that contains $G/U$.
If we allow enough minors to vanish, it's possible that no element of $G$ would satisfy the minor conditions.
However, we would need to take more than one frozen minor to zero for this to happen.
The frozen minors determine a normal crossing partial compactification of $G^{e,w_0}$, so this is a codimension 2 phenomenon.
That is, up to codimension 2, the canonical partial compactification of $G^{e,w_0}$ afforded by taking frozen variables to 0 is $G/U$.
These spaces support the same functions by Hartogs' Theorem.

The initial seed $s_0$ has cluster variables $A_{i;j}=\Delta^{1,\dots, i}_{j-i+1, \dots, j}$.
With the proof of Proposition~{\ref{prop:Gew0InGmodU}} in mind, 
this seems like a pretty natural choice of seed studying for functions on $G/U$.
As with Corollary~\ref{cor:W}, the potential $W$ for this seed can be written down immediately using Proposition~\ref{prop:theta}.

{\cor{\label{prop:W_Gew0}
Let $\cX$ be the Fock-Goncharov dual of the $\cA$-cluster variety $G^{e,w_0}$, 
and $W:\cX \to \C$ the potential arising from $G^{e,w_0}$'s frozen variables.
Then the restriction of $W$ to $T_{M;s_0}$ is
\eqn{
\sum_{i=1}^{n-1} \vartheta_{i;n} + \vartheta_{i;i} ,
} 
where 
\eqn{
\vartheta_{i;n} = \sum_{j=0}^{n-i-1} z^{-\sum_{k=0}^j e_{i;n-k}}
}
and
\eqn{
\vartheta_{i;i} = \sum_{j=0}^{i-1} z^{-\sum_{k=0}^j e_{i-k;i}}.
}
}}  

The basis of functions on $G/U$ corresponds to the integer points of the cone cut out by $W^T$.
A basis element $f$ of $H^0\lrp{G^{e,w_0}, \ssO_{G^{e,w_0}}}$ is also a basis element of $H^0\lrp{G/U, \ssO_{G/U}}$
if the $\gv$-vector of $f$ is a point in this cone. 
The integer points of the cone form a monoid whose non-identity generators are the primitive vectors along the edges of the cone.
We will see that the cone 
$\Xi := \lrc{x \in \cX\lrp{\R^T} : W^T(x) \geq 0  }$
is, upon a suitable change of basis, the Gelfand-Tsetlin cone $K_n$ with the final coordinate restricted to 0.
The edges of $K_n$ correspond to the minors $\Delta^{1,\dots,i}_{j_1,\dots,j_i}$.
On $\SL_n$, $\Delta^{1,\dots,n}_{1,\dots,n} \equiv 1$.
To compare directly to $K_n$ (instead of a slice of it),
we can add another frozen variable for this minor:  $A_{n;n} := \Delta^{1,\dots,n}_{1,\dots,n}$.
The double pseudoline arrangement yields no arrows to $v_{n;n}$,
so every seed is trivially optimized for $v_{n;n}$.
The new potential $\widetilde{W}$ is $W+ z^{-e_{n;n}}$.
$\widetilde{\Xi} := \lrc{x \in \cX\lrp{\R^T} : \widetilde{W}^T(x) \geq 0  }$
will be the Gelfand-Tsetlin cone.

{\remark{%
Note that I am not exactly performing the same analysis with $\GL_n$ instead of $\SL_n$.
While ${\GL_n}^{e,w_0}$ embeds into $\GL_n/U$,
its partial compactification 
$\overline{{\GL_n}^{e,w_0}}$ 
obtained by allowing all frozen variables to vanish 
does not agree with $\GL_n/U$ up to codimension 2.
The divisor $\Delta^{1,\dots,n}_{1,\dots,n} = 0$ has empty intersection with $\GL_n/U$.
The potential for $\GL_n$ would be the same as that for $\SL_n$,
but $\widetilde{W}$ is what actually yields the Gelfand-Tsetlin cone.
}}

Ideally, the edges of $\widetilde{\Xi}$ should correspond to $\gv$-vectors of the minors $\Delta^{1,\dots,i}_{j_1,\dots,j_i}$.
The $\gv$-vectors for the initial seed minors are known: $\gv_{s_0}\lrp{\Delta^{1,\dots,i}_{j-i+1,\dots,j}} = e_{i;j}^*$.
The remaining $\gv$-vectors can be computed iteratively.

{\prop{\label{prop:gvIteratively}
\begin{enumerate}
	\item $$\gv_{s_0} \lrp{ \Delta^{1,\dots,i+1}_{1,j-i+2,\dots,j+1} } = -e_{i;i+1}^* + e_{i;j+1}^* + e_{i+1;i+1}^*$$
	\item For $j_i < n$, if $$\gv_{s_0} \lrp{ \Delta^{1,\dots,i}_{j_1,\dots,j_i} }= \sum_{1\leq k \leq l \leq j_i} c_{k;l} \ e_{k;l}^*,$$
		then $$\gv_{s_0} \lrp{ \Delta^{1,\dots,i}_{j_1+1,\dots,j_i+1} }= \sum_{1\leq k \leq l \leq j_i} c_{k;l} \ e_{k;l+1}^*.$$
	\item For $j_i < n$, if $$\gv_{s_0} \lrp{ \Delta^{1,\dots,i}_{j_1,\dots,j_i} }= \sum_{1\leq k \leq l \leq j_i} c_{k;l} \ e_{k;l}^*,$$
		then $$\gv_{s_0} \lrp{ \Delta^{1,\dots,i+1}_{1,j_1+1,\dots,j_i+1} }= \sum_{1\leq k \leq l \leq j_i} c_{k;l} \ \gv_{s_0} 
		\lrp{\Delta^{1,\dots,k+1}_{1,l-k+2,\dots,l+1}}.$$
	\item Each $\gv_{s_0} \lrp{\Delta^{1,\dots,i}_{j_1,\dots,j_i}}$ can be computed using (1)-(3).
\end{enumerate}
}}

\begin{proof}
For {\it{(1)}}, it is convenient to use Proposition~\ref{prop:Gew0minors}.
Here $k=1$.
In (\ref{eq:mutijk}), the first summand of each binomial corresponds to outgoing arrows and the second to incoming arrows at $v_{i;j}$.
In this mutation sequence, at the vertex of mutation there are only outgoing arrows to $w_{i';j'}$ vertices.\footnote{%
Recall that $w_{i';j'}$ is the frozen vertex associated to the principal coefficient ${X_{i';j'}}_{s_0}$. 
The direction of these arrows is obtained from the proof of Lemma~\ref{lem:vWarrows}, {\it{mutatis mutandis}}.}
So only the first summand of each binomial obtains $X$ coefficients.
Over the central fiber, 
(\ref{eq:mutijk})
becomes
\eqn{
{A_{i;j}}_0\ \mu^*_{i;j}\lrp{{A_{i;j}}_1} = 
\begin{cases}
{A_{i;j+1}}_0\ {A_{i;j-1}}_1 & \mbox{if } j\neq i+1\\
{A_{i;j+1}}_0\ {A_{i+1;i+1}}_1 & \mbox{if } j= i+1.
\end{cases}
}
Then 
$\gv_{s_0} \lrp{\Delta^{1,\dots,i+1}_{1,3,\dots,i+2}} = -e_{i;i+1}^* + e_{i;i+2}^* + e_{i+1;i+1}^*$ and 
$\gv_{s_0} \lrp{\Delta^{1,\dots,i+1}_{1,j-i+1,\dots,j+1}} = \gv_{s_0} \lrp{\Delta^{1,\dots,i+1}_{1,j-i,\dots,j}}-e_{i;j}^* + e_{i;j+1}^*$.
Combining these yields {\it{(1)}}. 

For {\it{(2)}}, note that $\Delta^{1,\dots,i}_{j_1,\dots,j_i}$ is a cluster variable in ${\SL_{j_i}}^{e,w_0}$, 
so there is a mutation pattern producing this minor that uses only vertices $v_{k;l}$, $k<l< j_i$.
This mutation pattern describes $\Delta^{1,\dots,i}_{j_1,\dots,j_i}$ 
as a rational function of ${\SL_{j_i}}^{e,w_0}$'s initial seed minors.
Now note that replacing each vertex $v_{k;l}$ in the mutation pattern with $v_{k;l+1}$ 
simply shifts the column indices of every minor in the rational function up 1, 
{\it{i.e.}} $\Delta^{1,\dots,i}_{j-i+1,\dots,j}$ is replaced by
$\Delta^{1,\dots,i}_{j-i+2,\dots,j+1}$.
So the column indices of the resulting minor are shifted up 1 as well, 
as is the second index of each $e_{i';j'}^*$ in $\gv_{s_0}\lrp{\Delta^{1,\dots,i}_{j_1,\dots,j_i}}$.

I'll need a lemma for {\it{(3)}}.

{\lemma{\label{lem:minorMuPat}
Let the pair $\lrp{Q_\triangle,\mu_\triangle}$ denote the quiver and the sequence of mutations 
described in Lemma~\ref{lem:delbot},
with the number of rows of $Q_\triangle$ unspecified.
Note that the unfrozen subquiver of $G^{e,w_0}$'s initial seed has this form.
Within this quiver, for each $\Delta^{1,\dots,i}_{j_1,\dots,j_i}$ there is a descending chain of proper subquivers 
${Q_\triangle}_1\supset {Q_\triangle}_2 \supset \cdots \supset {Q_\triangle}_k$ 
such that the sequence of mutations ${\mu_\triangle}_1,{\mu_\triangle}_2,\dots, {\mu_\triangle}_k$
yields $\Delta^{1,\dots,i}_{j_1,\dots,j_i}$, the index corresponding to the bottom row of ${Q_\triangle}_1$ is
$j_i-1$,
and the index corresponding to the bottom row of ${Q_\triangle}_{r+1}$ is strictly less than 
the index corresponding to the bottom row of ${Q_\triangle}_{r}$ for all $r < k$.
}}

\begin{proof}
Choose ${Q_\triangle}_1$ such that the top vertex is $v_{1;j_1+1}$ and bottom right vertex is $v_{j_i-1-j_1;j_i-1}$.
After applying ${\mu_\triangle}_1$, the minor associated to the vertex $v_{k;l}$ of ${Q_\triangle}_1$ is
$\Delta^{1,\dots,k+1}_{j_1, l-k+2,\dots, l+1}$ by Proposition~\ref{prop:Gew0minors} together with {\it{(2)}}.
Recall from the proof of Proposition~\ref{prop:Gew0minors} that each mutation in the sequence is a direct application
of Lemma~\ref{lem:minorIdent}. 
Next take ${Q_\triangle}_2$ to have top vertex $v_{1;j_2}$ and bottom right vertex $v_{j_i-1-j_2;j_i-2}$.
If we had not already applied ${\mu_\triangle}_1$, the minor associated to $v_{k;l}$ in ${Q_\triangle}_2$ 
would have been 
$\Delta^{1,\dots,k+1}_{j_2-1, l-k+2,\dots, l+1}$.
But ${\mu_\triangle}_1$ has replaced every minor $\Delta^{1,\dots,k'}_{l'-k'+1,\dots,l'}$ 
used in ${\mu_\triangle}_2$  
with a new minor $\Delta^{1,\dots,k'+1}_{j_1,l'-k'+2,\dots,l'+1}$.
Note that shifting every column index by 1 and adding a new column ($j_1$) 
to every minor in Lemma~\ref{lem:minorIdent} preserves the equality.
So the minor associated to $v_{k;l}$ now must be  
$\Delta^{1,\dots,k+2}_{j_1,j_2,l-k+3,\dots, l+2}$.
We can continue this process.
The top vertex of ${Q_\triangle}_r$ will be $v_{1;j_r-r+2}$ and the bottom right vertex will be $v_{j_i-1-j_r;j_i-r}$.
The last valid ${Q_\triangle}_r$ occurs when $j_i-1-j_{r+1} <1$.
In this case, either $j_{r+1}=j_i$ or $j_{r+1}$ and $j_i$ are consecutive.
If $j_{r+1}=j_i$,
then $j_i-j_r \geq 2$ and 
$\Delta^{1,\dots,i}_{j_1,\dots, j_i}$
is obtained after ${\mu_\triangle}_r$'s mutation at 
$v_{1;j_i-r}$.
If $j_{r+1}+1=j_i$,
then $\Delta^{1,\dots,i}_{j_1,\dots, j_i}$
is obtained after ${\mu_\triangle}_r$'s mutation at 
$v_{2;j_i-r}$.
\end{proof}

So $\Delta^{1,\dots,i+1}_{1,j_1+1,\dots,j_i+1}$ can be obtained by a sequence of mutations $\mu$ that simply uses Lemma~\ref{lem:minorIdent} 
at each step.
Deleting the subsequence ${\mu_\triangle}_1$ from $\mu$ 
yields a sequence $\mu'$ that produces $\Delta^{1,\dots,i}_{j_1,\dots,j_i}$,
and the effect of ${\mu_\triangle}_1$ is to replace every minor $\Delta^{1,\dots,k}_{l-k+1,\dots,l}$ used in $\mu'$ with 
$\Delta^{1,\dots,k+1}_{1,l-k+2,\dots,l+1}$.
The column indices for an initial seed minor used in $\mu'$ cannot exceed $j_i$ as there are never arrows 
between the vertices used in $\mu'$ and vertices $v_{k;l}$ for $l>j_i$.
So $\gv_{s_0}\lrp{\Delta^{1,\dots,i}_{j_1,\dots,j_i}}$ has the form
\eq{
\gv_{s_0}\lrp{\Delta^{1,\dots,i}_{j_1,\dots,j_i}} &= \sum_{1\leq k \leq l \leq j_i} c_{k;l}\ 
\gv_{s_0}\lrp{\Delta^{1,\dots,k}_{l-k+1,\dots,l}}\\
&= \sum_{1\leq k \leq l \leq j_i} c_{k;l}\ e_{k;l}^*.
}{eq:gDelta}  
In both $\mu'$ and $\mu$, arrows between the vertex of mutation and some $w_{k;l}$ are always outgoing.
Then $\gv_{s_0}\lrp{\Delta^{1,\dots,i+1}_{1,j_1+1,\dots,j_i+1}}$ 
can be obtained by replacing each $\Delta^{1,\dots,k}_{l-k+1,\dots,l}$ in 
(\ref{eq:gDelta})
with $\Delta^{1,\dots,k+1}_{1,l-k+2,\dots,l+1}$,
yielding
\eqn{
\gv_{s_0}\lrp{\Delta^{1,\dots,i+1}_{1,j_1+1,\dots,j_i+1}} &= \sum_{1\leq k \leq l \leq j_i} c_{k;l}\ 
\gv_{s_0}\lrp{\Delta^{1,\dots,k+1}_{1,l-k+2,\dots,l+1}}.
}

For {\it{(4)}}, every increasing sequence $j_1, \dots, j_i$ of positive integers can be obtained 
from the sequences $l-k+1, \dots,l$ indexing the columns of the initial seed minors by some combination of the maps
$a: r_1, \dots, r_s \mapsto r_1 +1, \dots, r_s+1 $
and 
$b: r_1, \dots, r_s \mapsto 1, r_1 +1, \dots, r_s+1 $.
The $\gv$-vectors provided in {\it{(1)}} are precisely those needed to utilize {\it{(3)}}.
\end{proof}

The inequalities defining $\widetilde{\Xi}$ can be indexed by the frozen vertices and a subset of the arrows of the initial seed quiver.
I'll illustrate this for $n=4$.

\begin{minipage}{\linewidth}
\begin{center}
\begin{tikzpicture}
	\node (14) at (0,0) [draw, regular polygon, regular polygon sides=4, inner sep=0pt] {$v_{1;4}$};
	\node (24) at (2,0) [draw, regular polygon, regular polygon sides=4, inner sep=0pt] {$v_{2;4}$};
	\node (34) at (4,0) [draw, regular polygon, regular polygon sides=4, inner sep=0pt] {$v_{3;4}$};
	\node (44) at (6,0) [draw, regular polygon, regular polygon sides=4, inner sep=0pt] {$v_{4;4}$};

	\node (13) at (1,1) {$v_{1;3}$};
	\node (23) at (3,1) {$v_{2;3}$};
	\node (33) at (5,1) [draw, regular polygon, regular polygon sides=4, inner sep=0pt] {$v_{3;3}$};

	\node (12) at (2,2) {$v_{1;2}$};
	\node (22) at (4,2) [draw, regular polygon, regular polygon sides=4, inner sep=0pt] {$v_{2;2}$};

	\node (11) at (3,3) [draw, regular polygon, regular polygon sides=4, inner sep=0pt] {$v_{1;1}$};

	\draw[->] (14) -- (13) ;
	\draw[->] (13) -- (12) ;

	\draw[->] (24) -- (23) ;

	\draw[->] (33) -- (23) ;
	\draw[->] (23) -- (13) ;

	\draw[->] (22) -- (12) ;
\end{tikzpicture}
\captionof{figure}{Each box and each arrow represents an inequality defining $\widetilde{\Xi}$.
$\boxed{\protect\vphantom{{1^1_j}^1_j} v_{i;j}}$ corresponds to $\lra{e_{i;j},\ \cdot\ } \geq 0$.
$v_{i;j} \rightarrow v_{i;j-1}$  corresponds to $\lra{\sum_{k=j-1}^n e_{i;k},\ \cdot\ } \geq 0$, {\it{e.g.}} $v_{1;3} \rightarrow v_{1;2}$ corresponds to $\lra{e_{1;2} + e_{1;3} + e_{1;4}, \ \cdot \ } \geq 0$.
$v_{i;j} \leftarrow v_{i+1;j}$  corresponds to $\lra{\sum_{k=i}^j e_{k;j},\ \cdot\ } \geq 0$, {\it{e.g.}} $v_{1;3} \leftarrow v_{2;3}$ corresponds to $\lra{e_{1;3} + e_{2;3} + e_{3;3}, \ \cdot \ } \geq 0$.
\label{fig:IneqQuiver}}
\end{center}
\end{minipage}

Replacing the vertices of Figure~\ref{fig:IneqQuiver} with coefficients of a given $\gv$-vector 
provides a convenient way to see that all inequalities are satisfied and to see which are in fact equalities.
For example, for $n=4$, $\gv_{s_0}\lrp{\Delta^{1,2}_{1,3}} = -e_{1;2}^* + e_{1;3}^* + e_{2;2}^*$ would be represented by the picture below.

\input{IneqQuiverExample.tex}
 
An analogous picture is often drawn for the Gelfand-Tsetlin cone.
The inequalities in this case are represented by a triangular array like the one below.\cite{BKir,BZ96}

\input{GTIneq.tex}

$\Delta^{1,2}_{1,3}$ corresponds to the Gelfand-Tsetlin pattern having three $1$'s in the first downward diagonal, one $1$ in the second downward diagonal, and $0$'s elsewhere.

\input{GTIneqExample.tex}

Comparing pictures for several minors suggests a certain correspondence between the inequalities defining $\widetilde{\Xi}$ and those defining $K_n$.
I'll illustrate the correspondence in the $n=4$ case.

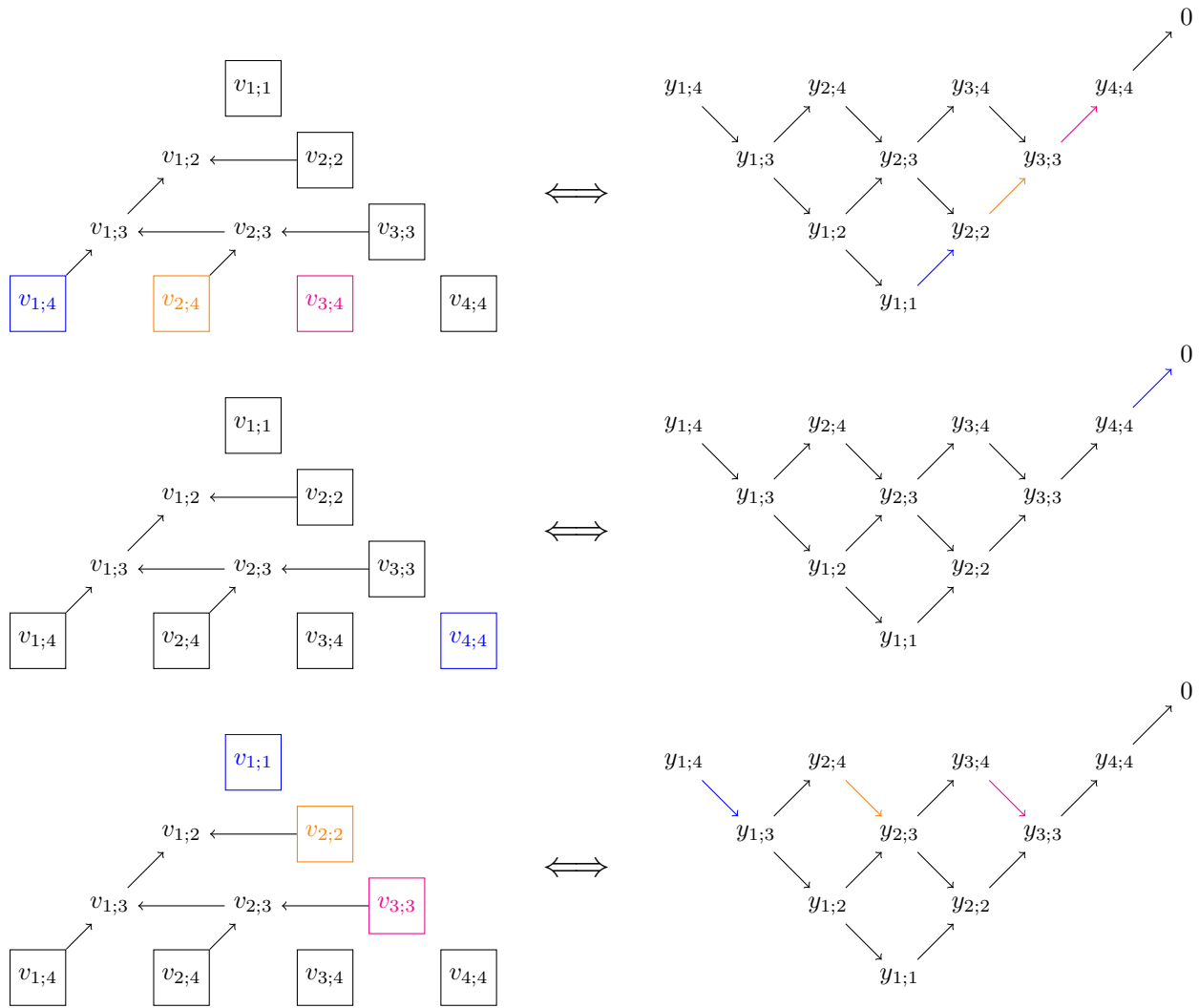
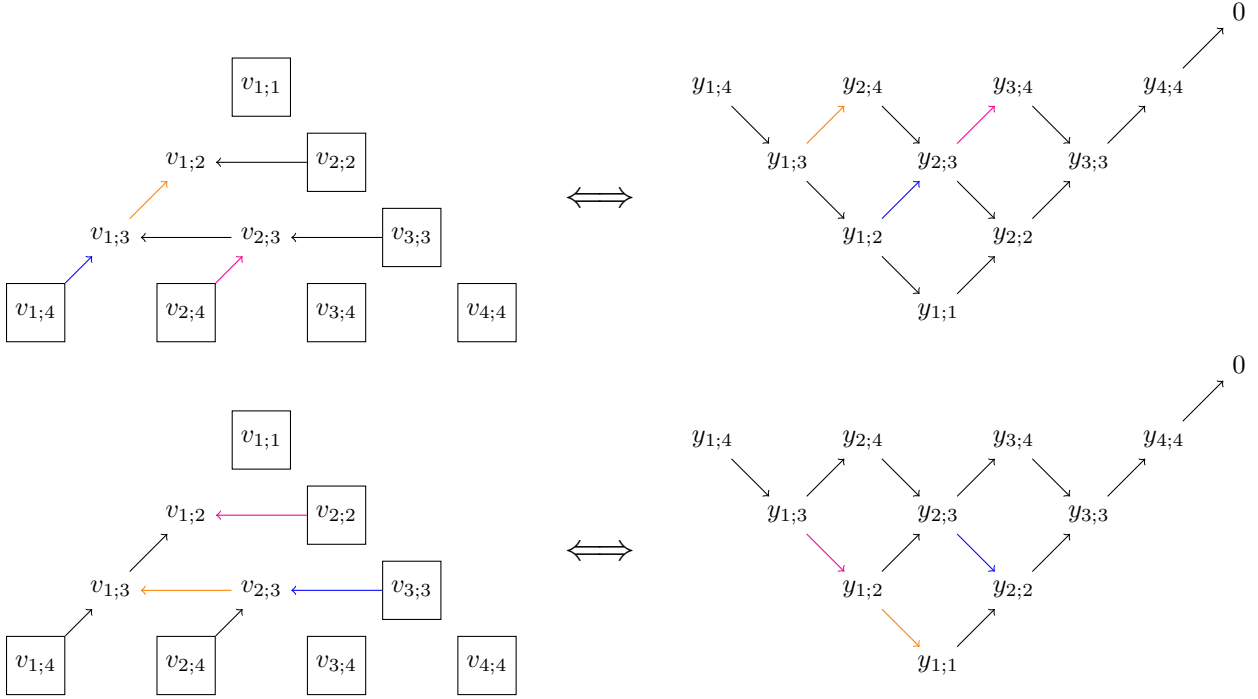
\input{IneqCorrespondence1.tex}
\input{IneqCorrespondence2.tex}
\input{IneqCorrespondence3.tex}
\input{IneqCorrespondence4.tex}
\input{IneqCorrespondence5.tex}

This apparent correspondence motivated the change of basis in the following theorem.

{\theorem{\label{thrm:XiGT}
Let 
$\psi: \R^{\frac{n(n+1)}{2}} \to \R^{\frac{n(n+1)}{2}}$ 
be the change of basis given by 
\eqn{
\psi \lrp{\sum_{1\leq i \leq j \leq n} y_{i;j}\ \vb{e}_{i;j}  } = \sum_{1\leq i \leq j \leq n} x_{i;j}\ \vb{e}_{i;j}',    
}
where
\eqn{
&x_{n;n} = y_{n;n}\\
&x_{i;i} = y_{i;n}- y_{i;n-1} \qquad {\text{for }} i<n\\
&x_{i;n} = y_{i;i}- y_{i+1;i+1} \qquad {\text{for }} i<n\\
&x_{i;j} = y_{i;n-j+i} -y_{i+1;n-j+i+1} + y_{i+1;n-j+i} - y_{i;n-j+i-1} \qquad {\text{for }} i<j<n.
}
Then $\psi\lrp{K_n} = \widetilde{\Xi}$.
%
}}

\begin{proof}
I need to verify that 
$\lrp{x_{1;1}, x_{1;2}, \dots, x_{n;n}}$ satisfies the inequalities defining $\widetilde{\Xi}$ 
precisely when 
$\lrp{y_{1;1}, y_{1;2}, \dots, y_{n;n}}$ satisfies the inequalities defining $K_n$.
Clearly, $x_{n;n} \geq 0 \iff y_{n;n} \geq 0$.
Next, for $i<n$, 
$x_{i;i} \geq 0 \iff y_{i;n} \geq y_{i;n-1}$
and 
$x_{i;n} \geq 0 \iff y_{i;i} \geq y_{i+1;i+1}$.
For $i<j<n$,
\eqn{
\sum_{k=j}^n x_{i;k} \geq 0 &\iff 
y_{i;i} - y_{i+1;i+1} + \sum_{k=j}^{n-1} \lrp{y_{i;n-k+i}- y_{i+1;n-k+i+1} + y_{i+1;n-k+i} - y_{i;n-k+i-1}}  
\geq 0\\
&\iff 
y_{i;i} - y_{i+1;i+1} + y_{i;n-j+i}- y_{i+1;n-j+i+1} + y_{i+1;n-(n-1)+i} - y_{i;n-(n-1)+i-1}  
\geq 0\\
&\iff 
y_{i;i} - y_{i+1;i+1} + y_{i;n-j+i}- y_{i+1;n-j+i+1} + y_{i+1;i+1} - y_{i;i}  
\geq 0\\
&\iff 
y_{i;n-j+i}- y_{i+1;n-j+i+1}
\geq 0\\
&\iff 
y_{i;n-j+i}\geq y_{i+1;n-j+i+1}
} 
and
\eqn{
\sum_{k=i}^j x_{k;j} \geq 0 &\iff 
y_{j;n} - y_{j;n-1} + \sum_{k=i}^{j-1} \lrp{y_{k;n-j+k}- y_{k+1;n-j+k+1} + y_{k+1;n-j+k} - y_{k;n-j+k-1}}  
\geq 0\\
&\iff 
y_{j;n} - y_{j;n-1} + y_{i;n-j+i}- y_{j;n} + y_{j;n-1} - y_{i;n-j+i-1}  
\geq 0\\
&\iff  
y_{i;n-j+i}  - y_{i;n-j+i-1}  
\geq 0\\
&\iff
y_{i;n-j+i} \geq  y_{i;n-j+i-1}.  
}
This accounts for all inequalities of both $\widetilde{\Xi}$ and $K_n$. $\psi\lrp{K_n} = \widetilde{\Xi}$.
\end{proof}

Since the minors $\Delta^{1,\dots,i}_{j_1,\dots,j_i}$ correspond bijectively to the edges of $K_n$, 
they must also be in bijection with the edges of $\widetilde{\Xi}$.
If their $\gv$-vectors are along distinct edges of $\widetilde{\Xi}$, 
then every edge of $\widetilde{\Xi}$ is $\R_{\geq 0}\cdot \gv_{s_0} \lrp{\Delta^{1,\dots,i}_{j_1,\dots,j_i}}$ 
for some $\Delta^{1,\dots,i}_{j_1,\dots,j_i}$.

{\prop{\label{prop:gvectEdges}
Let 
$GT\lrp{\Delta^{1,\dots,i}_{j_1,\dots,j_i}}$ 
denote the Gelfand-Tsetlin pattern associated to
$\Delta^{1,\dots,i}_{j_1,\dots,j_i} $.
Then 
\eqn{
\psi\lrp{GT\lrp{\Delta^{1,\dots,i}_{j_1,\dots,j_i}}} = \gv_{s_0}\lrp{\Delta^{1,\dots,i}_{j_1,\dots,j_i}}.
}}}

\begin{proof}
The $k;l$ component of 
$GT\lrp{\Delta^{1,\dots,i}_{j_1,\dots,j_i}}$
is 
\eqn{
GT\lrp{\Delta^{1,\dots,i}_{j_1,\dots,j_i}}_{k;l}= 
\begin{cases} 
1 & \mbox{if } j_{i+1-k} \geq n+1-l \\
0 & \mbox{otherwise.}
\end{cases}
}
Then for $k = n$, 
\eqn{
\psi\lrp{GT\lrp{\Delta^{1,\dots,i}_{j_1,\dots,j_i}}}_{n;n} &= 
\begin{cases} 
1 & \mbox{if } j_{i+1-n} \geq 1 \\
0 & \mbox{otherwise}
\end{cases}\\
&= 
\begin{cases} 
1 & \mbox{if } i = n\\
0 & \mbox{otherwise}.
\end{cases}
}
For $k<n$,
\eqn{
\psi\lrp{GT\lrp{\Delta^{1,\dots,i}_{j_1,\dots,j_i}}}_{k;k} &= 
GT\lrp{\Delta^{1,\dots,i}_{j_1,\dots,j_i}}_{k;n} -GT\lrp{\Delta^{1,\dots,i}_{j_1,\dots,j_i}}_{k;n-1}\\
&=
\begin{cases} 
1 & \mbox{if } j_{i+1-k} \geq 1\\ 
0 & \mbox{otherwise}
\end{cases}
-
\begin{cases} 
1 & \mbox{if } j_{i+1-k} \geq 2\\
0 & \mbox{otherwise}
\end{cases}\\
&= 
\begin{cases} 
1 & \mbox{if } j_{i+1-k} = 1\\
0 & \mbox{otherwise}
\end{cases}\\
&= 
\begin{cases} 
1 & \mbox{if } i=k \mbox{ and } j_{1} = 1\\
0 & \mbox{otherwise}
\end{cases}
}
and 
\eqn{
\psi\lrp{GT\lrp{\Delta^{1,\dots,i}_{j_1,\dots,j_i}}}_{k;n} &= 
GT\lrp{\Delta^{1,\dots,i}_{j_1,\dots,j_i}}_{k;k} -GT\lrp{\Delta^{1,\dots,i}_{j_1,\dots,j_i}}_{k+1;k+1}\\
&=
\begin{cases} 
1 & \mbox{if } j_{i+1-k} \geq n+1-k\\ 
0 & \mbox{otherwise}
\end{cases}
-
\begin{cases} 
1 & \mbox{if } j_{i-k} \geq n-k\\
0 & \mbox{otherwise}
\end{cases}\\
&=
\begin{cases} 
1 & \mbox{if } j_{i+1-k} = n+1-k\\ 
0 & \mbox{otherwise}
\end{cases}
-
\begin{cases} 
1 & \mbox{if } j_{i-k} = n-k\\
0 & \mbox{otherwise.}
\end{cases}
}
For $k<l<n$,
\eqn{
\psi\lrp{GT\lrp{\Delta^{1,\dots,i}_{j_1,\dots,j_i}}}_{k;l}&= 
GT\lrp{\Delta^{1,\dots,i}_{j_1,\dots,j_i}}_{k;n-l+k}
- 
GT\lrp{\Delta^{1,\dots,i}_{j_1,\dots,j_i}}_{k+1;n-l+k+1}\\
&\phantom{=}+
GT\lrp{\Delta^{1,\dots,i}_{j_1,\dots,j_i}}_{k+1;n-l+k}
-
GT\lrp{\Delta^{1,\dots,i}_{j_1,\dots,j_i}}_{k;n-l+k-1}\\
&=
\begin{cases} 
1 & \mbox{if } j_{i+1-k} \geq l+1-k \\
0 & \mbox{otherwise}
\end{cases}
-
\begin{cases} 
1 & \mbox{if } j_{i-k} \geq l-k \\
0 & \mbox{otherwise}
\end{cases}\\
&\phantom{=}+
\begin{cases} 
1 & \mbox{if } j_{i-k} \geq l+1-k \\
0 & \mbox{otherwise}
\end{cases}
-
\begin{cases} 
1 & \mbox{if } j_{i+1-k} \geq l+2-k \\
0 & \mbox{otherwise}
\end{cases}\\
&=
\begin{cases} 
1 & \mbox{if } j_{i+1-k} = l+1-k \\
0 & \mbox{otherwise}
\end{cases}
-
\begin{cases} 
1 & \mbox{if } j_{i-k} = l-k \\
0 & \mbox{otherwise.}
\end{cases}
}
Now consider the initial seed minors $\Delta^{1,\dots,i}_{j-i+1,\dots,j}$.
In this case, 
each of the expressions above reduces to 
\eqn{
\psi\lrp{GT\lrp{\Delta^{1,\dots,i}_{j-i+1,\dots,j}}}_{k;l} = 
\begin{cases} 
1 & \mbox{if } i = k \mbox{ and } j = l\\
0 & \mbox{otherwise,}
\end{cases}
}
in agreement with the $\gv$-vectors for these minors.
Next, consider the minors $\Delta^{1,\dots,i+1}_{1,j-i+2,\dots,j+1}$.
Now the expressions reduces to 
\eqn{
\psi\lrp{GT\lrp{\Delta^{1,\dots,i+1}_{1,j-i+2,\dots,j+1}}}_{k;k} = 
\begin{cases} 
1 & \mbox{if }  k = i+1\\
0 & \mbox{otherwise}
\end{cases}
}
and for $k<l$,
\eqn{
\psi\lrp{GT\lrp{\Delta^{1,\dots,i+1}_{1,j-i+2,\dots,j+1}}}_{k;l} = 
\begin{cases} 
1 & \mbox{if }  k = i,\ l = j+1, \mbox{ and } j \neq i\\
-1 & \mbox{if }  k = i,\ l =i+1, \mbox{ and } j \neq i\\
0 & \mbox{otherwise.}
\end{cases}
}
This also agrees with the $\gv$-vectors.

Next, I claim that for $j_i <n$,
$$\psi\lrp{GT\lrp{\Delta^{1,\dots,i}_{j_1+1,\dots,j_i+1}}}_{k;l+1}
=
\psi\lrp{GT\lrp{\Delta^{1,\dots,i}_{j_1,\dots,j_i}}}_{k;l}.$$
If $k<l$, this is clear.
If $k=l=n$, both are 0 since $i\leq j_i <n$.
If $k=l <n$,
$$\psi\lrp{GT\lrp{\Delta^{1,\dots,i}_{j_1,\dots,j_i}}}_{k;l}
=
\begin{cases}
1 & \mbox{if } i=k \mbox{ and } j_1 =1\\
0 & \mbox{otherwise,}
\end{cases}$$
while
\eqn{
\psi\lrp{GT\lrp{\Delta^{1,\dots,i}_{j_1+1,\dots,j_i+1}}}_{k;l+1}
&= 
\begin{cases}
1 & \mbox{if } j_{i+1-k} = 2\\
0 & \mbox{otherwise}
\end{cases}
-
\begin{cases}
1 & \mbox{if } j_{i-k} = 1\\
0 & \mbox{otherwise}
\end{cases}\\
&=
\begin{cases}
1 & \mbox{if } i=k \mbox{ and } j_1 = 1\\
0 & \mbox{otherwise.}
\end{cases}
}

Finally, I claim that for $j_i < n$,
\eqn{
\psi\lrp{GT\lrp{\Delta^{1,\dots, i+1}_{1,j_1+1, \dots, j_i+1}}}_{r;s}
=
\sum_{1\leq k \leq l \leq j_i} \psi\lrp{GT\lrp{\Delta^{1,\dots,i}_{j_1,\dots, j_i}}}_{k;l} \psi \lrp{GT\lrp{\Delta^{1,\dots, k+1}_{1,l-k+2, \dots, l+1}}}_{r;s}.
}
When $r=s$, we have
\eq{
\psi\lrp{GT\lrp{\Delta^{1,\dots, i+1}_{1,j_1+1, \dots, j_i+1}}}_{r;r}
=
\begin{cases}
1 & \mbox{if } r = i+1\\
0 & \mbox{otherwise,}
\end{cases}
}{eq:GTgenrr}
and 
\eqn{
\sum_{1\leq k \leq l \leq j_i} \psi\lrp{GT\lrp{\Delta^{1,\dots,i}_{j_1,\dots, j_i}}}_{k;l} \psi \lrp{GT\lrp{\Delta^{1,\dots, k+1}_{1,l-k+2, \dots, l+1}}}_{r;s}
&=
\sum_{1\leq k \leq l \leq j_i} \psi\lrp{GT\lrp{\Delta^{1,\dots,i}_{j_1,\dots, j_i}}}_{k;l} 
\cdot\begin{cases}
1 & \mbox{if } r = k+1\\
0 & \mbox{otherwise}
\end{cases}\\
&=
\sum_{r-1\leq l \leq j_i} \psi\lrp{GT\lrp{\Delta^{1,\dots,i}_{j_1,\dots, j_i}}}_{r-1;l}. 
}
This is 0 unless $r>1$, in agreement with (\ref{eq:GTgenrr}).  For $r>1$, 
\eqn{
\sum_{r-1\leq l \leq j_i} \psi\lrp{GT\lrp{\Delta^{1,\dots,i}_{j_1,\dots, j_i}}}_{r-1;l}
=&
\psi\lrp{GT\lrp{\Delta^{1,\dots,i}_{j_1,\dots, j_i}}}_{r-1;r-1}
+
\sum_{r-1 < l \leq j_i} \psi\lrp{GT\lrp{\Delta^{1,\dots,i}_{j_1,\dots, j_i}}}_{r-1;l}\\
=&
\begin{cases}
1 & \mbox{if } r = i+1 \mbox{ and } j_1 = 1\\
0 & \mbox{otherwise}
\end{cases}\\
&+
\sum_{r-1 < l \leq j_i} \lrp{
\begin{cases}
1 & \mbox{if } j_{i+2-r} = l+2-r\\
0 & \mbox{otherwise}
\end{cases}
-
\begin{cases}
1 & \mbox{if } j_{i+1-r} = l+1-r\\
0 & \mbox{otherwise}
\end{cases}
}\\
=& 
\begin{cases}
1 & \mbox{if } r = i+1 \mbox{ and } j_1 = 1\\
0 & \mbox{otherwise}
\end{cases}
+
\begin{cases}
1 & \mbox{if } r = i+1 \mbox{ and } j_1 \geq 2\\
0 & \mbox{otherwise}
\end{cases}\\
=&
\begin{cases}
1 & \mbox{if } r = i+1\\
0 & \mbox{otherwise.}
\end{cases}
}
This agrees with (\ref{eq:GTgenrr}) as well.
For $r<s$,
\eq{
\psi\lrp{GT\lrp{\Delta^{1,\dots, i+1}_{1,j_1+1, \dots, j_i+1}}}_{r;s}
=
\begin{cases}
1 & \mbox{if } j_{i+1-r}  = s-r\\
0 & \mbox{otherwise}
\end{cases}
-
\begin{cases}
1 & \mbox{if } j_{i-r} +1  = s-r \mbox{ or } \lrp{ r = i \mbox{ and } s = r+1} \\
0 & \mbox{otherwise}
\end{cases}
}{eq:GTgenrs}
and
\eq{
&\sum_{1\leq k \leq l \leq j_i} \psi\lrp{GT\lrp{\Delta^{1,\dots,i}_{j_1,\dots, j_i}}}_{k;l} \psi \lrp{GT\lrp{\Delta^{1,\dots, k+1}_{1,l-k+2, \dots, l+1}}}_{r;s}\\
&=
\sum_{1\leq k \leq l \leq j_i} \psi\lrp{GT\lrp{\Delta^{1,\dots,i}_{j_1,\dots, j_i}}}_{k;l} 
\cdot
\begin{cases}
1 & \mbox{if } r = k,\ s= l+1, \mbox{ and } k \neq l \\
-1 & \mbox{if } r = k,\ s= k+1, \mbox{ and } k \neq l \\
0 & \mbox{otherwise}
\end{cases}\\
&=
\sum_{ r < l \leq j_i} \psi\lrp{GT\lrp{\Delta^{1,\dots,i}_{j_1,\dots, j_i}}}_{r;l} 
\cdot
\begin{cases}
1 & \mbox{if }  s= l+1 \\
-1 & \mbox{if } s= r+1 \\
0 & \mbox{otherwise}
\end{cases}\\
&=
\begin{cases}
\psi\lrp{GT\lrp{\Delta^{1,\dots,i}_{j_1,\dots, j_i}}}_{r;s-1} & \mbox{if }  s= l+1 \mbox{ for some } l,\ r< l \leq j_i  \\
-\psi\lrp{GT\lrp{\Delta^{1,\dots,i}_{j_1,\dots, j_i}}}_{r;r} & \mbox{if }  s= r+1 \\
0 & \mbox{otherwise}
\end{cases}\\
&=
\begin{cases}
\begin{cases}
1 & j_{i+1-r} = s -r\\ 
0 & \mbox{otherwise}
\end{cases}
-
\begin{cases}
1 & j_{i-r} = s -r -1 \\
0 & \mbox{otherwise}
\end{cases} & \mbox{if } r+1 < s \leq j_i\\
-1 & \mbox{if } r= i,\ j_1 = 1, \mbox{ and } s = r+1 \\
0 & \mbox{otherwise.}
\end{cases}
}{eq:GTgenrsprime}
(\ref{eq:GTgenrs}) and (\ref{eq:GTgenrsprime}) agree as well.
Then by Proposition~\ref{prop:gvIteratively},
\eqn{
\psi\lrp{GT\lrp{\Delta^{1,\dots,i}_{j_1,\dots, j_i}}} = \gv_{s_0} \lrp{\Delta^{1,\dots,i}_{j_1,\dots, j_i}} 
} 
for each $\lrp{\Delta^{1,\dots,i}_{j_1,\dots, j_i}}$. 
\end{proof}

{\cor{\label{cor:gvectsPrimGens}}
The $\R_{\geq 0}$ span of the $\gv$-vectors $\gv_{s_0} \lrp{\Delta^{1,\dots,i}_{j_1,\dots,j_i}}$ form the edges of $\widetilde{\Xi}$.
Furthermore, these $\gv$-vectors are primitive in $\widetilde{\Xi} \cap \cX\lrp{\Z^T} \subset \cX\lrp{\R^T}$.
}

\begin{proof}
This follows immediately from Proposition~\ref{prop:gvectEdges}
since the analogous statements for the Gelfand-Tsetlin cone are known, 
and a subset of these $\gv$-vectors form a basis for the $\cX\lrp{\Z^T}$ (identified with a lattice by the choice of seed).
\end{proof}

As discussed in \cite{GHKK}, Proposition~\ref{prop:opt_seeds}, Proposition~\ref{prop:gvectEdges}, and Corollary~\ref{cor:gvectsPrimGens} together imply the full Fock-Goncharov conjecture for $G/U$. 
Recall that the open embedding $G^{e,w_0} \hookrightarrow G/U$
used in this section was given by $g \mapsto g^T U$.
So the minor $\Delta^{1,\dots, i}_{j_1,\dots, j_i}$ on $G^{e,w_0}$ 
is identified with the minor $\Delta^{j_1,\dots, j_i}_{1,\dots,i}$ on $G/U$.  	  
Then the canonical basis $B_{G/U}$ of $H^0\lrp{G/U, \ssO_{G/U}}$
is parametrized by $\Xi \cap \cX\lrp{\Z^T}$, and its elements can be expressed (not necessarilly uniquely) as
\eqn{
\prod_{\lrc{j_1 < \cdots < j_i} \subsetneq \lrc{1,\dots,n}} \lrp{ \Delta_{1,\dots,i}^{j_1,\dots,j_i} }^{a_{j_1,\dots,j_i}},
}
where
$ a_{j_1,\dots,j_i} \in \Z_{\geq 0}$.

The maximal torus $H$ in $\SL_n$ determined by the choice of opposite Borel subgroups $B_+$ and $B_-$ is the subgroup of diagonal matrices.
The torus action on $H^0\lrp{G/U, \ssO_{G/U}}$ is given by $h \cdot f\lrp{g U} = f\lrp{g h U}$.
Each $f\in B_{G/U}$ is an eigenfunction of this action. 
Each $h\in H$ has the form $\diag\lrp{h_1, \dots, h_n}$, $h_1 h_2 \cdots h_n = 1$.
Right multiplication by a diagonal matrix $h$ scales column $j$ by $h_j$.
So $h\cdot \Delta^{j_1,\dots,j_i}_{1,\dots,i} = h_1 h_2\cdots h_i\ \Delta^{j_1, \dots, j_i}_{1,\dots,i}$. 
We can record the weight of an eigenfunction of this action by recording the exponents of each $h_i$.
So $\Delta^{j_1, j_2, j_3}_{1,2,3}$ has weight $\lrp{1,1,1,0,\dots, 0}$, 
and $\Delta^{k}_{1} \Delta^{j_1,j_2}_{1,2}$ has weight $\lrp{2,1,0,\cdots,0}$.
Fixing a weight $\lambda= \lrp{\lambda_1,\lambda_2, \dots, \lambda_{n-1}}$
intersects $\Xi$ with an affine hyperplane $\mathcal{H}_{\lambda_i}$ for each $i$.
The integer points of 
\eqn{\Xi\bigcap\lrp{\cap_{i=1}^{n-1} \mathcal{H}_{\lambda_i}}}
form a basis for the weight space $V_\lambda$,
and the number of integer points in this slice is then the dimension of $V_\lambda$.
The next proposition provides a sanity check, verifying that the desired weights are being produced.

{\prop{If $\gv_{s_0} \lrp{\Delta^{j_1,\dots,j_i}_{1,\cdots,i}} = \sum c_{k;l}\ e_{k;l}^*$,
then 
\eqn{
\sum_{l=k}^n c_{k;l} = 
\begin{cases}
1 & \mbox{if } i = k\\
0 & \mbox{otherwise.}
\end{cases} 
} 
}}

\begin{proof}
The coefficients $c_{k;l}$ are described in Proposition~\ref{prop:gvectEdges}.
For $k>i$, neither $j_{i+1-k}$ nor $j_{i-k}$ are valid row indices.
So $c_{k;l} = 0$ for all $l$ when $k > i$.
For $k<i$,\footnote{Note that each expression in Proposition~\ref{prop:gvectEdges} can be written in this form, possibly by adding 0.}   
\eqn{
c_{k;l} = 
\begin{cases}
1 & \mbox{if } j_{i+1-k} = l+1 - k\\ 
0 & \mbox{otherwise}
\end{cases}
-
\begin{cases}
1 & \mbox{if } j_{i-k} = l - k\\ 
0 & \mbox{otherwise.}
\end{cases}
}
Since $k<i$, both $j_{i+1-k}$ and $j_{i-k}$ exist.
Some $l$ with $k\leq l \leq n$ must satisfy $j_{i+1-k} = l+1-k$,
and this $l$ is obviously unique.
Similarly, some $l'$ (perhaps $l'=l$) with $k\leq l' \leq n $ must satisfy $j_{i-k} = l-k$,
and this $l'$ is unique.
Then $l$ contributes $+1$ to the sum and $l'$ contributes $-1$.
$\sum_{l=k}^n c_{k;l} = 0$ for $k<i$.

Finally, if $k=i$
\eqn{
c_{k;l} = 
\begin{cases}
1 & \mbox{if } j_{1} = l+1 - i\\ 
0 & \mbox{otherwise.}
\end{cases}
}
$l = j_1 +i -1$ contributes $1$ to the sum, while all other terms contribute 0.
\end{proof}

{\example{Let $n=3$, and consider the weight $\lambda=\lrp{3,1}$.
Then a point $x= \sum_{i\leq j } x_{i;j} e_{i;j}^*$ in $\cX\lrp{\R^T}$
is in $\mathcal{H}_{\lambda_1}\cap \mathcal{H}_{\lambda_2}$ 
if and only if
\eqn{
x_{2;2}+ x_{2;3} = 1
}
and
\eqn{
x_{1;1} + x_{1;2} + x_{1;3} = 3-1 = 2.
}
$x$ is in $\Xi\cap \mathcal{H}_{\lambda_1} \cap \mathcal{H}_{\lambda_2}$ 
if additionally
$x_{1;1} \geq 0$, 
$x_{2;2} \geq 0$, 
$x_{2;2} + x_{1;2} \geq 0$, 
$x_{1;3} \geq 0$, 
$x_{1;3} + x_{1;2} \geq 0$, 
and
$x_{2;3} \geq 0$.
The only integer points satisfying this system of equalities and inequalies are the ones listed below.

\begin{center}
\begin{tabular}{l | c c c c c}
{ } &	$x_{1;1}$ & $x_{1;2}$ & $x_{1;3}$ & $x_{2;2}$ & $x_{2;3}$ \\ \hline
1   &	0 	& -1	  & 3	    & 1	      & 0       \\
2   &	0 	& 0 	  & 2	    & 0	      & 1       \\
3   &	0 	& 0 	  & 2	    & 1	      & 0       \\
4   &	0 	& 1	  & 1	    & 0	      & 1       \\
5   &	0 	& 1	  & 1	    & 1	      & 0       \\
6   &	0 	& 2	  & 0	    & 0	      & 1       \\
7   &	0 	& 2	  & 0	    & 1	      & 0       \\
8   &	1 	& -1	  & 2	    & 1	      & 0       \\
9   &	1 	& 0	  & 1	    & 0	      & 1       \\
10  &	1 	& 0	  & 1	    & 1	      & 0       \\
11  &	1 	& 1	  & 0	    & 0	      & 1       \\
12  &	1 	& 1	  & 0	    & 1	      & 0       \\
13  &	2 	& -1	  & 1	    & 1	      & 0       \\
14  &	2 	& 0	  & 0	    & 0	      & 1       \\
15  &	2 	& 0	  & 0	    & 1	      & 0       \\
\end{tabular}
\end{center}

So $V_\lambda$ is a 15 dimensional irreducible representation of $\SL_3$.
}}

\nocite{*}
\bibliography{bibliography}
\bibliographystyle{hep}

\end{document}